\documentclass[a4paper,10pt]{article}
\usepackage[utf8]{inputenc}

\usepackage{authblk}

\usepackage[margin=2cm]{geometry}

\usepackage{xcolor}

\usepackage{graphicx}

\usepackage{float}

\usepackage{longtable}
\usepackage{ctable}% For making \specialrule in tabular or supertabular environments (see Sec:mathematical model, Table: Summary of model.
\usepackage{color, colortbl}
\usepackage{tabularx}
\newcolumntype{L}[1]{>{\raggedright\arraybackslash}p{#1}}
\newcolumntype{C}[1]{>{\centering\arraybackslash}p{#1}}
\newcolumntype{R}[1]{>{\raggedleft\arraybackslash}p{#1}}
\newcounter{Eqnno}
\DeclareRobustCommand{\Eqn}[1]{%
   \refstepcounter{Eqnno}%
   \theEqnno\label{#1}}

\title{Testing a thermo-chemo-hydro-geomechanical model for gas hydrate bearing sediments using triaxial compression lab experiments}
\author[1]{S. Gupta}
\author[1]{C. Deusner}
\author[1]{M. Haeckel}
\author[2]{R. Helmig}
\author[3]{B. Wohlmuth}

\affil[1]{
Helmholtz Centre for Ocean Research Kiel,
Wischhofstr. 1-3, D-24148 Kiel, Germany.
}

\affil[2]{
Dept. of Hydromechanics and Modelling of Hydrosystems,
University of Stuttgart,
Pfaffenwaldring 61, 
70569 Stuttgart,
Germany.
}

\affil[3]{
Chair for Numerical Mathematics, 
Technical University Munich,
Boltzmannstr. 3, 
85748 Garching bei M\"unchen, 
Germany.
}
	  
\date{}

\begin{document}

\maketitle

\begin{abstract}
Natural gas hydrates are considered a potential resource for gas production on industrial scales. Gas hydrates contribute to the strength and stiffness of the hydrate-bearing sediments. During gas production, the geomechanical stability of the sediment is compromised. Due to the potential geotechnical risks and process management issues, the mechanical behavior of the gas hydrate-bearing sediments needs to be carefully considered.
In this study, we describe a coupling concept that simplifies the mathematical description of the complex interactions occuring during gas production by isolating the effects of sediment deformation and hydrate phase changes. Central to this coupling concept is the assumption that the soil grains form the load-bearing solid skeleton, while the gas hydrate enhances the mechanical properties of this skeleton. We focus on testing this coupling concept in capturing the overall impact of geomechanics on gas production behavior though numerical simulation of a high-pressure isotropic compression experiment combined with methane hydrate formation and dissociation.  We consider a linear-elastic stress-strain relationship because it is uniquely defined and easy to calibrate. Since, in reality, the geomechanical response of the hydrate bearing sediment is typically inelastic and is characterized by a significant shear-volumetric coupling, we control the experiment very carefully in order to keep the sample deformations small and well within the assumptions of poro-elasticity. The closely co-ordinated experimental and numerical procedures enable us to validate the proposed simplified geomechanics-to-flow coupling, and set an important precursor towards enhancing our coupled hydro-geomechanical hydrate reservoir simulator with more suitable elasto-plastic constitutive models.
\end{abstract}

\paragraph{Key points}
\begin{enumerate}
\item A simplified coupling concept is proposed for modelling geomechanical feedback on gas production behaviour of hydrate reservoirs.
\item Coupling concept is validated by numerical simulation of gas hydrate formation and dissociation in a controlled triaxial compression test. 
\item Coupling concept applies to pore-filling hydrates as well as for hydrates in the transition zone between pore-filling and load-bearing habits. 
\end{enumerate}

%% ------------------------------------------------------------------------ %%
%
%  TEXT
%
%% ------------------------------------------------------------------------ %%

\section{Introduction}
\label{sec:introduction}
Methane hydrates are crystalline solids formed from water molecules enclathrating methane molecules. 
Methane hydrates are thermodynamically stable under conditions of low temperatures and high pressures. 
If warmed or depressurized, methane hydrates destabilize and dissociate into water and methane gas. 
Natural gas hydrates occur in permafrost regions and the deep sea, usually in soils or sediments at considerable depth when methane is available in sufficient amounts. 
Natural gas hydrates are considered to be a promising energy resource. 
It is widely believed that the energy content of methane occurring in hydrate form is immense, possibly exceeding the combined energy content of all other conventional fossil fuels (\cite{Pinero2013}, \cite{Burwicz2011}, \cite{Archer2009}, \cite{Milkov2004}, \cite{Kvenvolden1993}).

Several methods have been proposed for production of natural gas from gas hydrate reservoirs, e.g., thermal stimulation, depressurization, and chemical activation (\cite{Moridis2009,Moridis2011}, \cite{Park2006}, \cite{Lee2003}). 
Currently, depressurization is deemed the most mature approach. Consequently, significant research and development effort has been directed towards assessing the potential of depressurization as a primary driving force for natural gas production from gas hydrate reservoirs. 
Recent field trials, onshore below the Alaskan permafrost and in the Nankai Trough offshore Japan were both essentially depressurization tests; 
the Japanese test used only depressurization (\cite{Yamamoto2013,Yamamoto2015}, \cite{David2013}), 
while, the Alaskan test was combined with $N_2$:$CO_2$ injection (\cite{Anderson2014}, \cite{Schoderbek2013}).

In the earlier gas production studies, several mathematical models (e.g. \cite{Tsypkin1991}, \cite{Ahmadi2004}, \cite{Yousif1991}, \cite{SunMohanty2005}, \cite{LiuFlemming2007}, \cite{Moridis2003,Moridis2007}) 
and numerical simulators (e.g. MH21-HYDRES \cite{MH21-HYDRES}, STOMP-HYD \cite{STOMP-HYD}, UMSICHT HyRes \cite{UMSICHT-HyRes}, TOUGH-HYDRATE \cite{TOUGH-HYDRATEmanual}) were developed, which focused on hydrate phase change and fluid flow rather than on the geomechanical behaviour. 
Over the years, it has become increasingly clear that the geomechanical effects associated with these gas production methods cannot be ignored. 
% The gas hydrate dissociation and fluid flow can lead to serious geomechanical stability issues like soil consolidation, seafloor subsidence and wellbore collapse, and in extreme cases, can potentially trigger catastrophic events like underwater landslides (\cite{SultanCochonat2004a,SultanCochonat2004b}).
Recent field trials have shown that large deformation and sand production are relevant risks for natural gas production from gas hydrate-bearing sediments (\cite{Schoderbek2013}, \cite{Yamamoto2013,Yamamoto2015}), and reliable simulation tools are needed for risk assessment and production strategy development. 
Coupling between solid deformation and fluid transport lays the foundation for the simulation of the thermo-hydro-chemo-mechanical behavior of gas hydrate-bearing sediment during gas production, and the experimental validation of the coupling relationships is extremely important for adding certainty to predictive simulation of production scenarios and sediment mechanical behavior in general.
% It is widely accepted in the gas hydrate community that to make any meaningful evaluation of future gas production technologies from gas hydrate reservoirs, it is necessary to include detailed risk quantification from the inherent geohazards. 
Several mathematical and numerical tools (e.g., \cite{KlarSogaNG2010}, \cite{Kimoto2010}, \cite{Rutqvist2011}, \cite{Hyodo2014}, \cite{SGupta2015}) have since been developed to study gas production in gas hydrate reservoirs in a coupled thermo-chemo-hydro-geomechanical framework.

\begin{figure}
 \centering
 \includegraphics[scale=0.55]{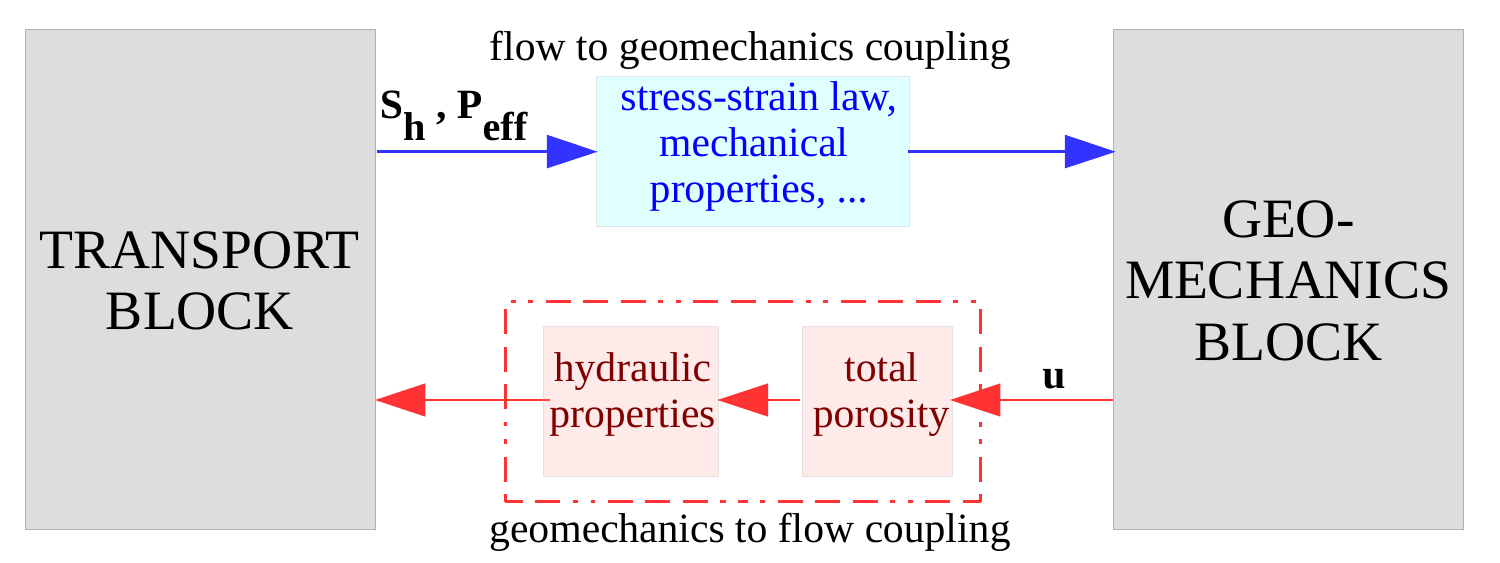}
 \caption{Coupling concept}
 \label{fig:coupling-concept}
\end{figure}

In a typical gas hydrate reservoir, the structure of the sediment is expected to change due to two distinct effects: 
1) the changing hydrate saturation, and 2) the sediment deformation. 
What complicates the matter further is that the hydrate provides additional strength to the sediment through a cementation-like effect, thereby, effectively coupling the two inputs; hydrate saturation, and sediment deformation. 
For any detailed hydro-geomechanical description of the gas production from gas hydrate-bearing sediments, it is imperative to analyze how the transport processes (i.e., flow and chemical processes) would respond to any given geomechanical input. 
As can be expected, this is a rather challenging task due to the complexity of the interactions. 
In our model, we simplify the mathematical description of the coupled hydro-geomechanical processes by conceptualizing that the model can be decomposed into two distinct model blocks: transport-block and geomechanics-block, with the coupling between the two manifesting as changes in properties of each model block. (See Fig. \ref{fig:coupling-concept}.) 
The transport-block solves for the hydrate phase change and the non-isothermal, two-phase, two-component flow of water and methane gas, while the geomechanics-block solves for the sediment displacements. 
This decomposition is based on the simplifying assumption that the soil grains constitute the skeleton of the porous matrix, while the gas hydrate enhances the mechanical properties of this skeleton without actively bearing the load. 
The relative deformation of the gas hydrate phase with respect to the soil skeleton is ignored. 
This assumption allows us to distinguish between the \textit{total porosity} and the \textit{apparent porosity}. 
The total porosity characterizes the total pore volume, i.e. the volume not occupied by the soil grains, 
while the apparent porosity characterizes the actual pore volume which is available for the fluid flow. 
The deformation of the hydrate-bearing sediment directly affects only the total porosity. 
The evolution of the actual or apparent porosity field is then modelled by scaling the total porosity with functions of hydrate saturation through simple geometric arguments. 
To make the physical meaning clear, the apparent porosity is the actual measured quantity, while the total porosity is a mathematical construct which allows us to isolate the effects of sediment deformation from those of hydrate phase change. 
We also assume that those properties of the transport-block which depend on the sediment structure (i.e., hydraulic properties like permeability, capillary pressure, specific surface area, tortuosity etc.) can be modelled as a \textit{multiplicative decomposition} of functions of total porosity and hydrate saturation. 
In effect, with these simplifications, we can describe all feedback from geomechanics-block to the transport-block through a single transfer variable: total porosity, and this forms a central feature of our coupling concept.

In this study, through numerical simulation of a highly controlled high-pressure triaxial experiment combined with methane hydrate formation and dissociation, we aim to establish whether this coupling concept is effective in capturing the overall impact of geomechanics on the gas production behaviour.

To be able to test the coupling with confidence, there are two important pre-requisites: 
1) a well-tested model for the transport-block, and 2) a good enough estimation of the displacement field. 
The validation of the transport-block in our hydrate reservoir model was performed in our earlier study (\cite{SGupta2015}). 
For the second pre-requisite, however, we require a suitable constitutive model for describing the stress-strain response of the hydrate-bearing sediment.
A number of non-linear elastic (e.g. \cite{YuSong2011nonlinear,Miyazaki2011nonlinear,Miyazaki2012nonlinear}), elasto-plastic (e.g. \cite{KlarSogaNG2010,UchidaSogaYamamoto2012MHCS,SunGao2015,Lin2015SMP}), and elasto- viscoplastic (e.g. \cite{Kimoto2010}) constitutive models have been proposed in the recent years to model the geomechanical behaviour of gas hydrate-bearing sediments. 
% The stress strain relationships in these models form highly non-linear, non-smooth, and dynamic optimization problems for which uniqueness and existence of an optimal solution can be guaranteed only within the assumptions of small-strains \textcolor{black}{with strong constraints on work-hardening (See \citep{lubliner2013plasticity}, Sec. 3.4)}. 
\textcolor{black}{The stress strain relation in these models is quite complex. 
Variationally consistent formulations result in non-linear and non-smooth inequality settings.
These can be reformulated in terms of non-linear complementarity functions to which semi-smooth Netwon algorithms as iterative solvers can be applied, which converge locally. 
To enlarge the local convergence radius, suitable regularisation and damping strategies can be designed (\cite{WohlmuthHager2009,WohlmuthHager2010}).
Only in special situations existence and uniqueness is given, and often, only local existence is guaranteed and path dependent solutions exist.
We refer to \cite{Alexander2004,Alexander2009} and the references therein for existence results at finite strain. 
However, none of these results can be directly applied to our setting since we have a fully coupled hydrate system involving more non-linearities in the bi-directional couplings.}
Additionally, the constraints of the system involve inequalities which are nonlinear and non-smooth, leading to a model with a large number of parameters which often makes model calibration challenging and unreliable. 
To the best of our knowledge, none of these models have been validated in coupled hydro-geomechanical settings in the context of gas production from gas hydrate reservoirs, and have a large uncertainty associated with their predictive capabilities in highly dynamic conditions. 
In this study, it is of particular interest to reduce the complexity of the geomechanics-block as much as possible in order to reduce the uncertainty associated with the choice of a constitutive model. 
We, therefore, chose a uniquely invertible linear-elastic constitutive model in the geomechanics-block of our hydrate reservoir model. 
We account for the stiffening effect due to gas hydrates by parameterizing the Young’s modulus as a function of hydrate saturation $S_h$ (\cite{SantamarinaRuppel2010}). 
We also account for material compressibility with respect to hydrostatic pressure. 
The pressure dependence of compressibility and the $S_h$ dependence of Lame’s parameters introduces a weak nonlinearity in the geomechanics-block. 
The numerical implementation of the poro-elastic model and the transport-to-geomechanics coupling in our hydrate reservoir model was also validated in our earlier study (\cite{SGupta2015}).

In general, poro-elasticity is not a realistic model for the geomechanical description of cemented granular materials where the stress-strain response is typically non-linear, and the shear-volumetric coupling (dilatancy) is of particular importance.
Therefore, in order to validate our coupling concept within the constraints stated above, we control our triaxial experiment very carefully in a way that ensures that the assumptions of poroelasticity remain valid throughout the periods of interest for numerical simulation.

\section{Hydrate reservoir model}
\label{sec:hydrateReservoirSimulator}

\begin{figure}
  \centering
    \includegraphics[scale=0.38]{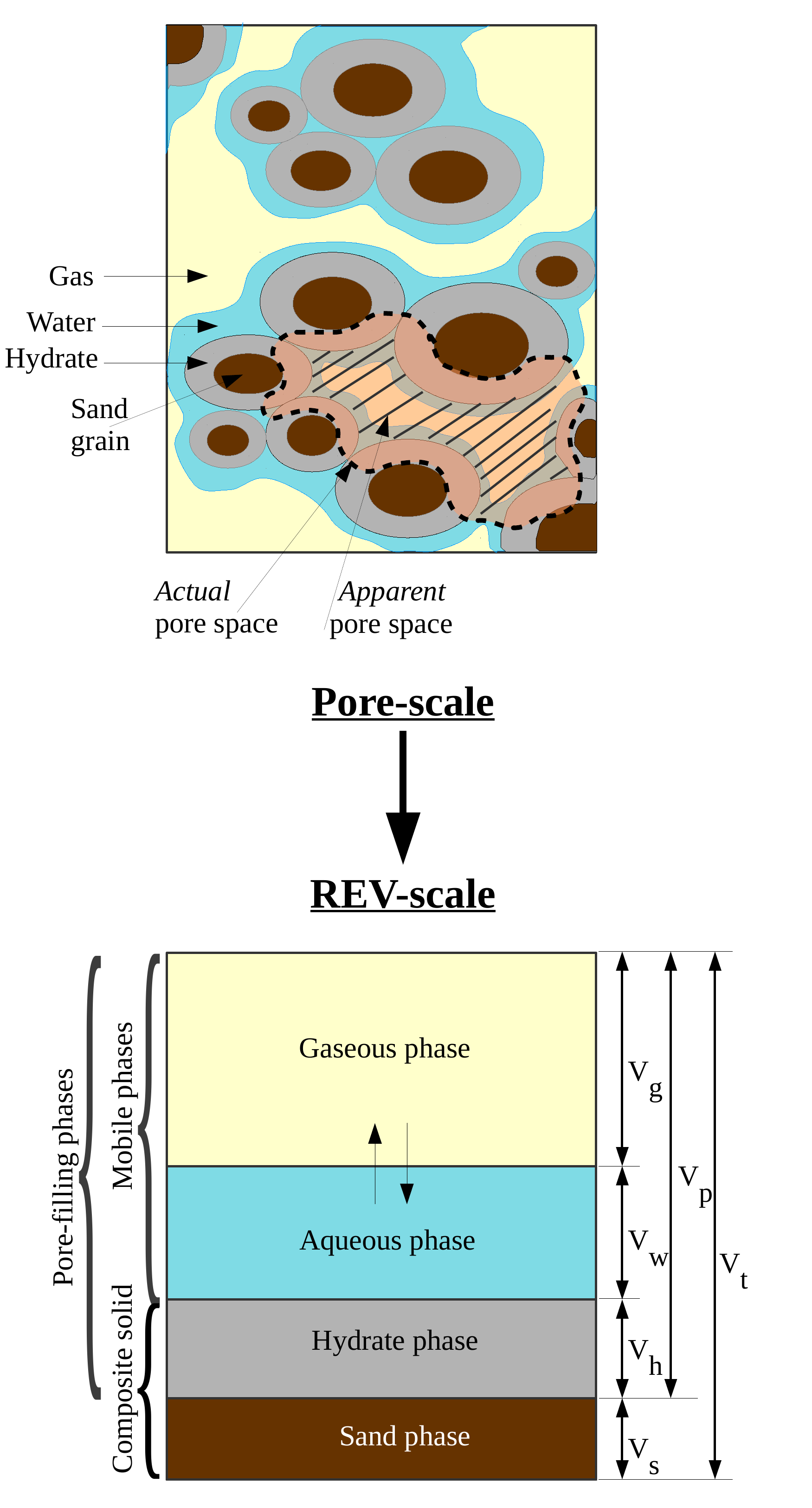}
  \caption{Pore-scale to REV-scale}  \label{fig:REVDescription}
\end{figure}

\paragraph{Mathematical model}
Our model considers kinetic hydrate phase change and non-isothermal, multi-phase, multi-component fluid flow through poro-elastic porous medium.
The model assumes that the porous medium is composed of three \textit{components}: $CH_4$, $H_2O$ and methane hydrate ($CH_4\cdot{N_h} H_2O$), 
which are present in three distinct \textit{phases}: gaseous, aqueous, and solid.
The gaseous phase comprises of methane gas and water vapour. 
The aqueous phase comprises of water and dissolved methane.
The solid phase comprises of pure methane hydrate and sand grains.
The sand grains are assumed to form a material continuum which provides the skeletal structure to the porous medium. We shall refer to this as \textit{solid matrix}.
The aqueous, gaseous, and hydrate phases exist in the void spaces of this solid matrix (See Fig.\ref{fig:REVDescription}).

We assume that the hydrate cements the sand grains in the mechanical sense (i.e., without forming any chemical bonds), such that the sand and hydrate together form a \textit{composite solid matrix}. The relative deformation of the gas hydrate phase with respect to the soil skeleton is ignored.
This assumption allows us to write a single momentum balance equation for the composite solid phase, instead of separate ones for the sand and the hydrate phases each.
To describe the mechanical behaviour of the composite solid matrix, we make a further simplifying assumption that the sand grains form the primary load-bearing structure, while the hydrates enhance the mechanical strength and stiffness of this structure without bearing any load themselves.
\textcolor{black}{This assumption is also adopted even after hydrate-bearing soil is loaded during depressurization. 
Because of this, the constitutive law does not consider stress-relaxation term that accounts for the release of load that has been carried by hydrates as introduced by others (\cite{KlarSogaNG2010,UchidaSogaYamamoto2012MHCS}).}

\textcolor{black}{
We decompose the mathematical model into transport and geomechanics blocks, and isolate the effects of hydrate phase change and ground deformation through the introduction of the variable \textit{total porosity}. 
This is justified based on the above assumption that the soil grains constitute the primary load-bearing skeleton of the porous matrix, such that deformation of the hydrate-bearing sediment directly affects only the total porosity.
This assumption allows us to solve for the mass balance of soil and hydrates separately, thus conveniently separating hydrate phase change kinetics from sediment deformation. 
If, for example, we do not make this assumption, then we would have to solve for the mass balance of soil and hydrate as a single composite phase leading to a strong coupling between hydrate phase change and geomechanics. 
The model decomposition would not be straightforward in this case, and the evolution of porosity field would be very complex, making the description of the fluid flow and the evolution of the hydraulic properties also significantly more challenging. 
One clear advantage of this simplification is that it gives a very modular structure to the model, with each model-block operating independently, and communicating with each other through coupling relationships which are neatly resolved with respect to the independent output variables of each model-block. (See Fig. 1.) 
Another important advantage is that the model decomposition allows us to use multi-rate time-stepping schemes, as discussed in \cite{SGupta-MRT2016}, which can significantly speed up the calculation, especially for 3D problems.
}

The mathematical model is described in detail in \cite{SGupta2015}.
A summary of the governing equations and the constitutive relationships is given in Table \ref{table:modelSummary}. %Table \ref{table:governingEquations} and Table \ref{table:constitutiveRelationships}.
The phases occupying the pore space (gaseous, aqueous, hydrate) are denoted by '$\beta $' $ = g, w, h $ respectively, 
the mobile phases (gaseous and aqueous) are denoted by '$\alpha$' $ = g, w $, 
and the mobile molecular components are denoted by '$\kappa$' $ = CH_4, H_2O $. 
The solid matrix is designated with the subscript '$s$'. 
The sand+hydrate composite solid matrix is designated with the subscript '$sh$'.
'$\gamma$' is used to denote all phases, i.e., $\gamma = g,w,h,$ and $s$.

\paragraph{Primary variables}
The mathematical model consists of the following $6$ governing equations: 
the mass balance Eqns. (\ref{eqn:MassBal_g},\ref{eqn:MassBal_w},\ref{eqn:MassBal_h},\ref{eqn:MassBal_s}), the momentum balance Eqn. (\ref{eqn:MomentumBal_sh}) for the composite-solid, and the energy balance Eqn. (\ref{eqn:EnergyBal}).
The momentum balance Eqns. (\ref{eqn:Darcy_g},\ref{eqn:Darcy_w}) for the mobile phases $\alpha=g,w$ give the $\alpha-$phase velocities directly, and are thus absorbed in the mass and energy balance equations.
We chose the following set of variables as the primary variables: 
the gas phase pressure $P_g$, the aqueous phase saturation $S_w$, the hydrate phase saturation $S_h$, the temperature $T$, the total-porosity $\phi$, and, the composite-matrix displacement $\mathbf{u}$.
All other variables can be derived (explicitly or implicitly) from this set of variables using the closure and constitutive relationships.

\paragraph{Solution strategy}
We use an iteratively coupled solution strategy.
The mathematical model is decomposed into three parts:
\begin{enumerate}
 \item transport-block ($\mathbf{F_{f}}$), comprised of the mass balance equations for $CH_4$, $H_2O$, and Hydrate, and the energy balance equation,   
 \item geomechanical-block ($\mathbf{F_{g}}$), comprised of the momentum balance equation for composite solid phase, and
 \item porosity-equation ($\mathbf{F_{\phi}}$), comprised of the mass balance equation for the sand phase.
\end{enumerate}

$\mathbf{F_{f}}$ is solved for the variables $P_g$, $S_w$, $S_h$, and $T$, $\mathbf{F_{g}}$ is solved for displacements $\mathbf{u}$, and $\mathbf{F_{\phi}}$ is solved for total porosity $\phi$.
$\mathbf{F_{f}}$ and $\mathbf{F_{\phi}}$ are spatially discretized using a fully upwinded cell centered finite volume method.
Orthogonal grids aligned with the principal axes are defined and a control-volume formulation with two-point flux approximation (TPFA) is used.
$\mathbf{F_{g}}$ is discretized using Galerkin finite element (FEM) method defined on $Q1$ elements.
An implicit Euler time-stepping scheme is used for marching forward in time.
The solution for a given time step involves two iterative loops, the \textit{inner loop} and the \textit{outer loop}.
The \textit{inner loop} uses Newton's method and SuperLU (\cite{SuperLU}) linear solver to solve each of $\mathbf{F_{f}}$, $\mathbf{F_{g}}$, and $\mathbf{F_{\phi}}$, thus taking care of the decoupled solution.
The \textit{outer loop} re-introduces the coupling between $\mathbf{F_{f}}$, $\mathbf{F_{g}}$, and $\mathbf{F_{\phi}}$ through a block Gauss Seidel iterative scheme.

The numerical scheme is implemented in the C++ based DUNE-PDELab framework (\cite{DUNEPdelab}), and is capable of solving problems in $1D$, $2D$ and $3D$ domains.

\section{Material and Methods}
\label{sec:experimentDetails}
We performed a controlled triaxial volumetric strain test on a sand sample in which methane hydrate was first formed under controlled effective stress and then dissociated via depressurization under controlled total stress. 
Gas hydrate in our experiment was initially formed by pressurizing partially water-saturated sand with gaseous methane to reach a gas hydrate saturation of $0.4$, and remaining methane gas was replaced with seawater before the sample was depressurized stepwise. 
Confining and axial loads in the triaxial testing were applied isotropically and were carefully controlled to keep the deformation of the sample small and well within the assumptions of poroelasticity.

\begin{figure}
  \centering
    \parbox{\textwidth}{a) \textcolor{black}{Gas hydrate formation.}}
    \vfill
    \includegraphics[scale=0.55]{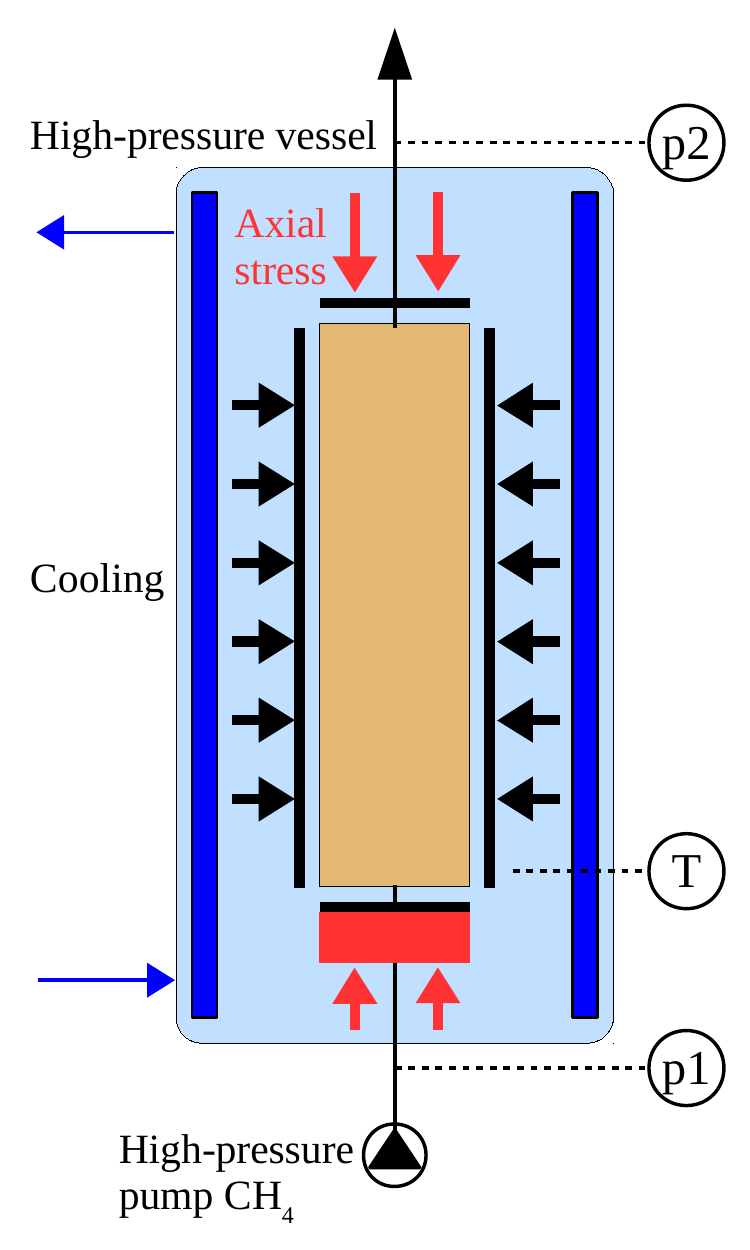}
    \vfill
    \parbox{\textwidth}{b) \textcolor{black}{Depressurization and gas production.}}
    \vfill
    \includegraphics[scale=0.55]{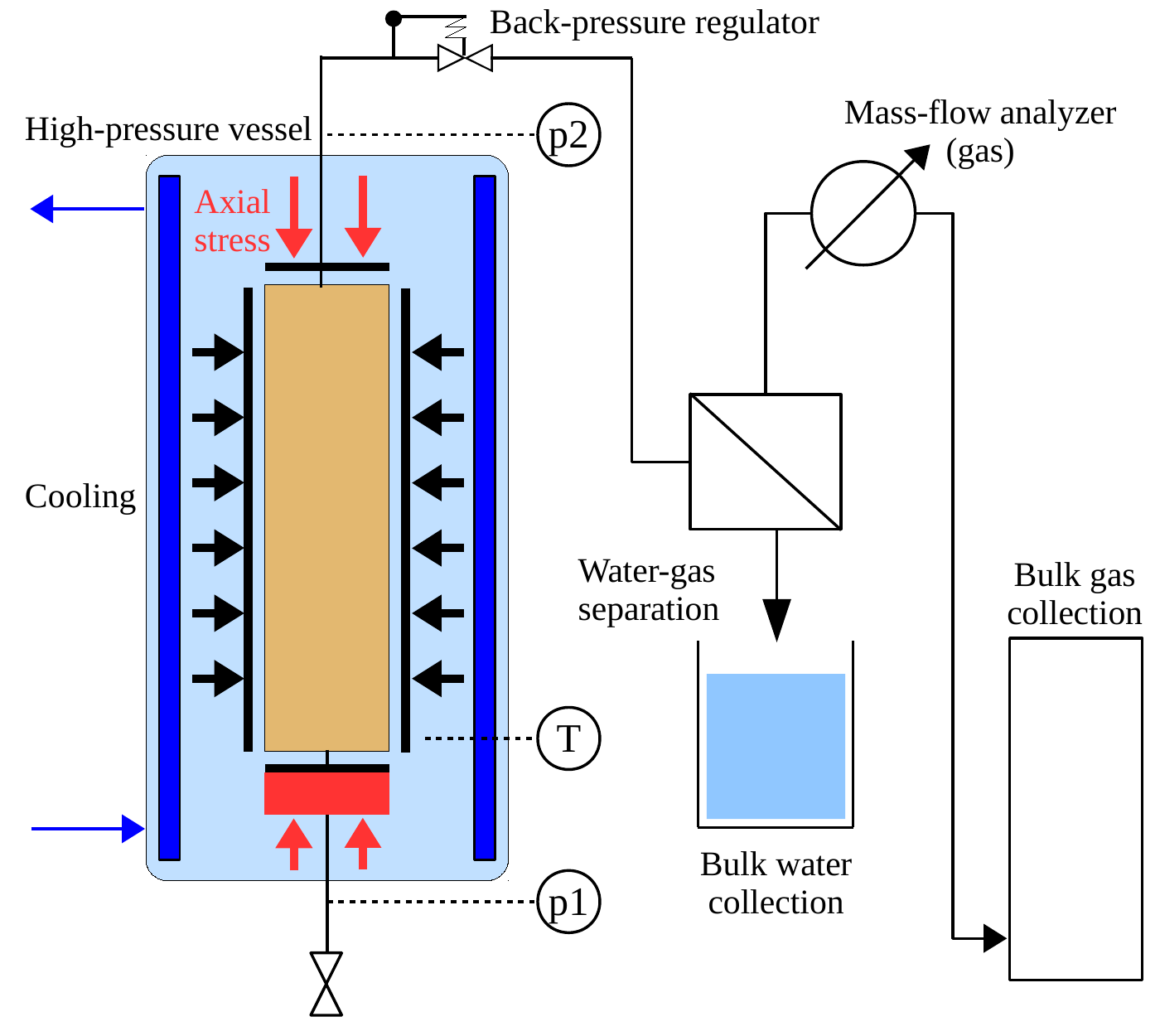}
    \caption{Simplified flow schemes for relevant period. }
    \label{fig:experimentalsetup}
\end{figure}

\paragraph{Experimental set-up and components} 
Experiments were carried out in the custom-made high-pressure apparatus NESSI (Natural Environment Simulator for Sub-seafloor Interactions, \cite{Deusner2012}), 
which is equipped with a triaxial cell mounted in a $40$ L stainless steel vessel (APS GmbH Wille Geotechnik, Rosdorf, Germany). 
The sample sleeve is made from FKM, other wetted parts of the setup are made of stainless steel. 
Salt water medium was stored in reservoir bottles (DURAN, Wertheim, Germany) prior to use in experiments, and the seawater medium was pressurized in an additional pressure vessel (Parr Instrument GmbH, Frankfurt, Germany) to allow fast transfer into the sample vessel. 
Fluid pressure in the sample vessel was adjusted with a back-pressure regulator valve (TESCOM Europe, Selmsdorf, Germany). 
Experiments were carried out in upflow mode with injection of $CH_4$ gas and seawater medium at the bottom of the sample prior to and after gas hydrate formation (Fig.\ref{fig:experimentalsetup}-a), respectively, and fluid discharge at the top of the sample during depressurization (Fig.\ref{fig:experimentalsetup}-b). 
\textcolor{black}{Axial and confining stresses, and the sample volume changes were monitored throughout the experimental period. 
Axial and confining stresses were controlled using high-precision hydraulic pumps and actuators (VPC $400$, APS GmbH Wille Geotechnik, Rosdorf, Germany), 
and changes in hydraulic fluid volumes were converted to calculate sample volume changes.}
Pore pressure was measured in the influent and the effluent fluid streams close to the sample top and bottom. 
The experiment was carried out under constant temperature conditions. 
Temperature control was achieved with a thermostat system (T1200, Lauda, Lauda-K\"onigshofen, Germany). 
Produced gas mass flow was analyzed with mass flow controllers (EL FLOW, Bronkhorst, Kamen, Germany). 
For control purposes, bulk effluent fluids were also collected inside $100$ L gas tight TEDLAR sampling bags (CEL Scientific, Santa Fe Springs CA, USA). 
The sampling bags were mounted inside water filled containers. 
After expansion of the effluent fluids at atmospheric pressure, the overall volume was measured as volume of water displaced from these containers.

\paragraph{Sample preparation and mounting}
The sediment sample was prepared from quartz sand (initial sample porosity: $0.35$, grain size: $0.1 - 0.6$ mm, G20TEAS, Schlingmeier, Schw\"ulper, Germany),
which was mixed with de-ionized water to achieve a final water saturation of $0.4$ relative to the initial sample porosity. 
The partially water-saturated and thoroughly homogenized sediment was filled into the triaxial sample cell to obtain final sample dimensions of $360$ mm in height and $80$ mm in diameter. 
The sample geometry was assured using a sample forming device. 
The sample was cooled to $2^0$C after the triaxial cell was mounted inside the pressure vessel. 
Initial water permeability of the gas hydrate-free sediment was estimated to be $5\times 10^{-10}$ $m^2$ .

\paragraph{Experimental procedure}

\paragraph{a) Gas hydrate formation}
\textcolor{black}{Prior to the gas hydrate formation, the partially water-saturated sediment sample was isotropically consolidated to $1$ MPa effective stress under drained conditions. 
It should be noted, that the apparent effective stress is monitored and controlled as differential pressure between the confining hydraulic fluid pressure of the reactor and the gas pressure in the sample pore space. 
Measurements and control algorithms do not take into account the changes in effective stress due to changes in water saturation and capillary pressure.
Thus, only the apparent effective stress is directly accessible from experimental procedures.}
The sample was flushed with $CH_4$ gas \textcolor{black}{to replace air with methane.} 
\textcolor{black}{The sample was}, subsequently, pressurized with $CH_4$ gas to approximately $12.5$ MPa (Fig.\ref{fig:experimentalsetup}-a). 
During pressurization with $CH_4$ gas and throughout the overall gas hydrate formation period, formation stress condition of $1$ MPa effective stress were maintained using an
automated control algorithm. 
The formation process was continuously monitored by logging the $CH_4$ gas pressure.
Mass balances and volume saturations were calculated based on $CH_4$ gas pressure and initial mass and volume values. 
Gas hydrate formation was terminated after $1.84$ mol of $CH_4$-hydrate had been formed after approximately $6$ days, corresponding to CH4-hydrate saturation of $0.39$. 
The sample was cooled to $-5^0$C and stress control was switched to constant total isotropic stress control at approximately $9$ MPa before the sample pore space was de-pressurized to atmospheric pressure and the remaining $CH_4$ gas in the pore space was released. 
System re-pressurization and water saturation of the pore space was achieved by instant filling and re-pressurization with pre-cooled $-1^0$C saltwater medium according to the seawater composition. 
Hydrate dissociation during the brief period of depressurization was minimized by taking advantage of the anomalous self-preservation effect, which reaches an optimum close to the chosen temperature (\cite{Stern2003}). 
After completion of the gas-water fluid exchange, the sample temperature was re-adjusted to $2^0$C. 

\paragraph{b) Depressurization and gas production}
The sample pore space was de-pressurized and gas produced by stepwise decrease of back pressure at constant isotropic total stress (Fig.\ref{fig:experimentalsetup}-b).
Overall fluid production (water and $CH_4$ gas) was monitored after de-pressurization at \textcolor{black}{atmospheric} pressure after temperature equilibration.

\section{Numerical Simulation}
\label{sec:numericalSimulation}
\begin{figure}
 \centering
 \includegraphics[scale=0.6]{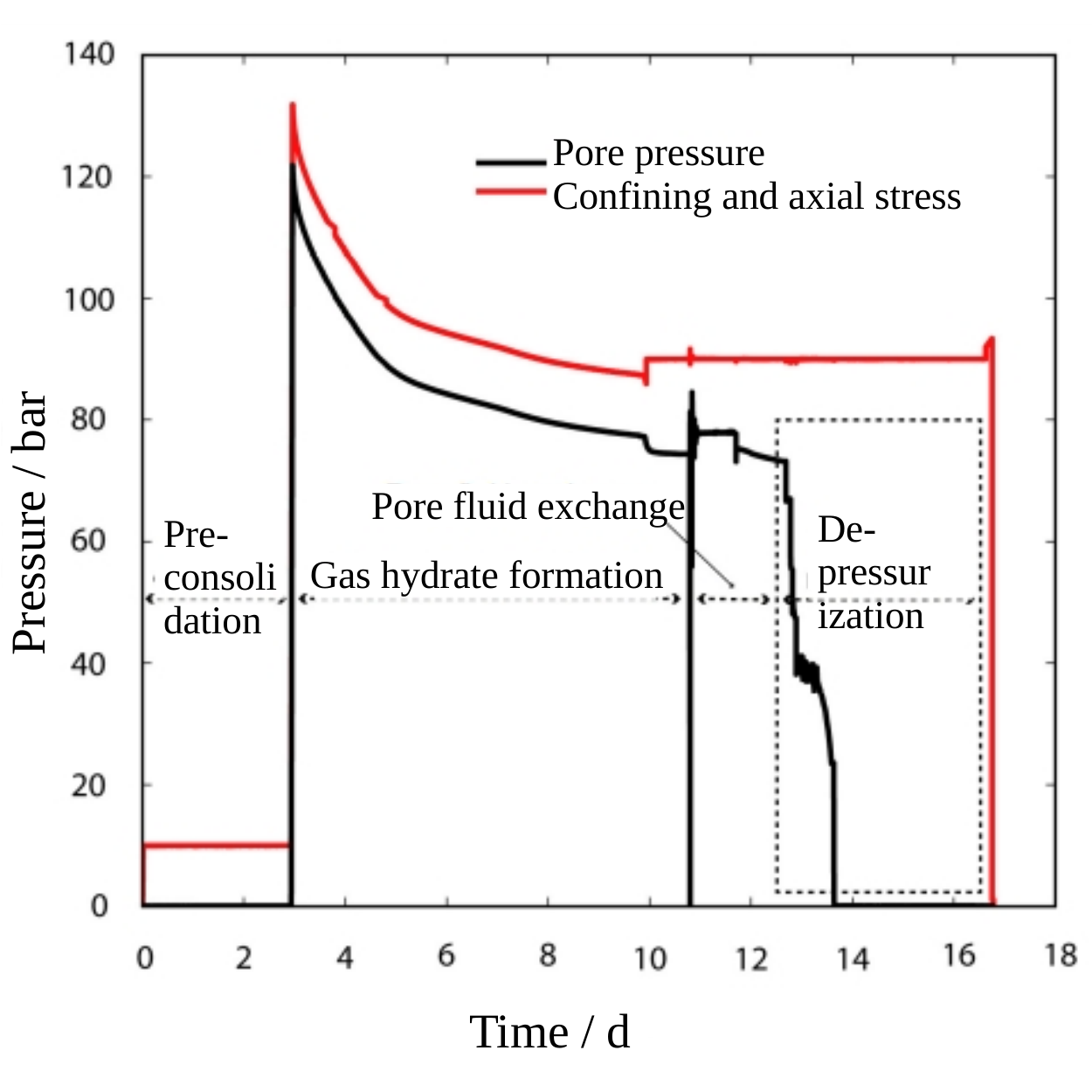}
 \caption{Overview of the \textit{measured} pressure and stresses over time.}
  \label{fig:pressureOverview}
\end{figure}

The overall experiment was carried out in four steps, viz. 1) \textit{pre-consolidation}, 2) \textit{gas hydrate formation}, 3) \textit{pore-fluid exchange}, and 4) \textit{depressurization}, as described in section \ref{sec:experimentDetails}. 
During steps $1$ and $2$, the sample  was maintained under a defined effective loading with the confining and the axial stresses were controlled to remain $10$ bar above the pore pressure. 
During steps $3$ and $4$, the total isotropic stress was controlled to remain at a constant level.  (See Fig. \ref{fig:pressureOverview}.)
The experiment was performed over a total period of about $16.8$ days. 
The \textit{periods of interest} for this simulation are: 1) from $Day-3$ to $Day-10$, corresponding to gas hydrate formation, and 2) from $Day-12.8$ to $Day-13.8$, corresponding to depressurization and gas production.
We simulate both of these periods separately.

% \label{fig:pressureOverview}

\paragraph{Computational domain}
Assuming that the sand sample is axially symmetric, a $2D$ radial plane of dimensions $360$ mm $\times$ $40$ mm is chosen as the computational domain. The dimensions correspond to the physical size of the sample.
The domain is discretized into $72 \times 8$ cells.

\begin{figure*}
 \centering
 \includegraphics[scale=0.65]{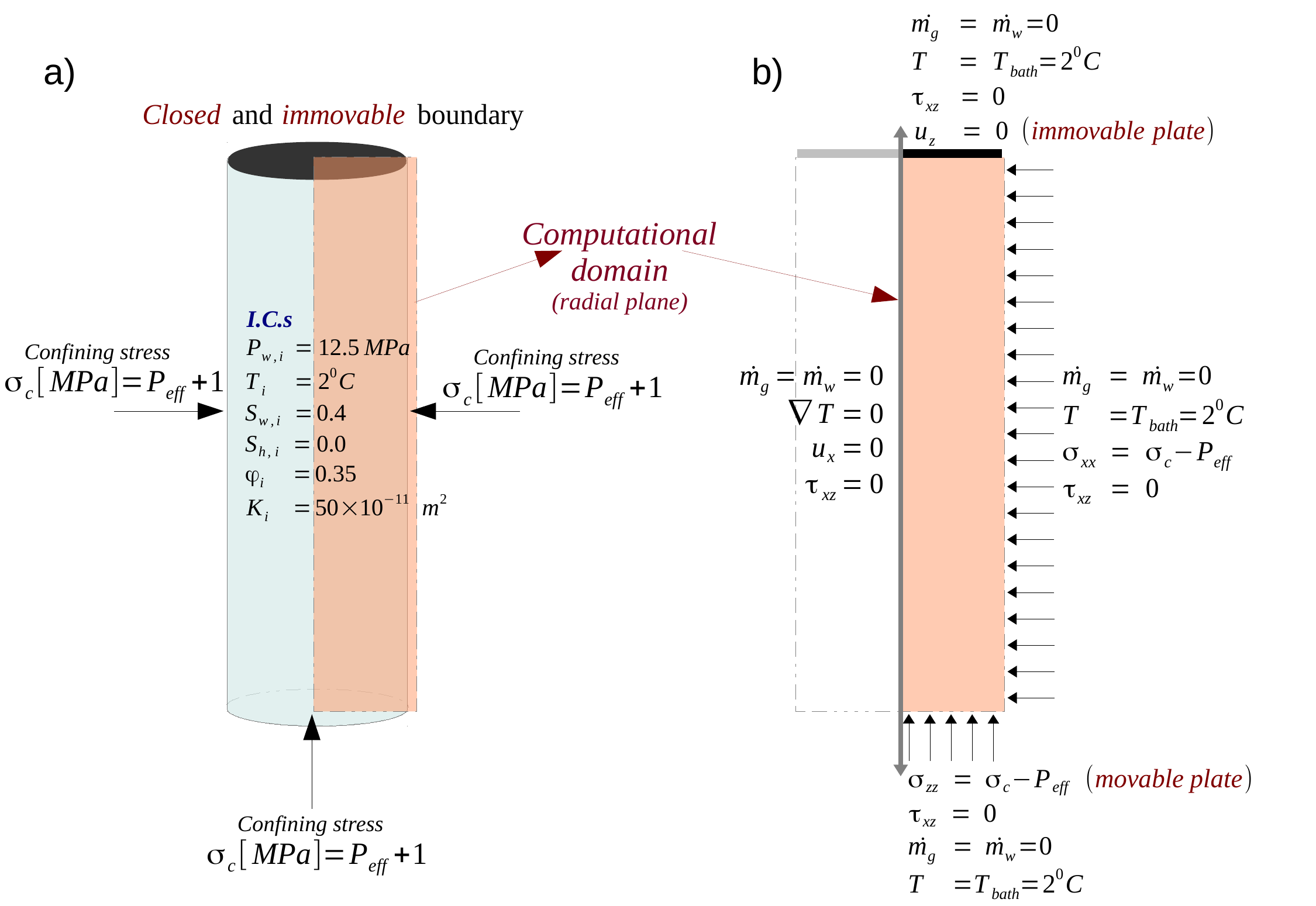}
 \caption{Test setting for the \textbf{gas hydrate formation} period. \newline
	  a) shows the sample and the initial conditions, and b) shows the $2D$ computational domain and the boundary conditions.}
  \label{fig:formation_testSetting}
\end{figure*}

\paragraph{Test-setting}
\paragraph{a) Gas hydrate formation period}
The schematic of the hydrate formation test is shown in Fig.\ref{fig:formation_testSetting}.
The schematic also shows the initial and boundary conditions.
The simulation is run until $t_{end}=604800$ s (i.e. $7$ days) using a maximum time step size of $120$ s.

% \label{fig:formation_testSetting}

\begin{figure*}
 \centering
 \includegraphics[scale=0.65]{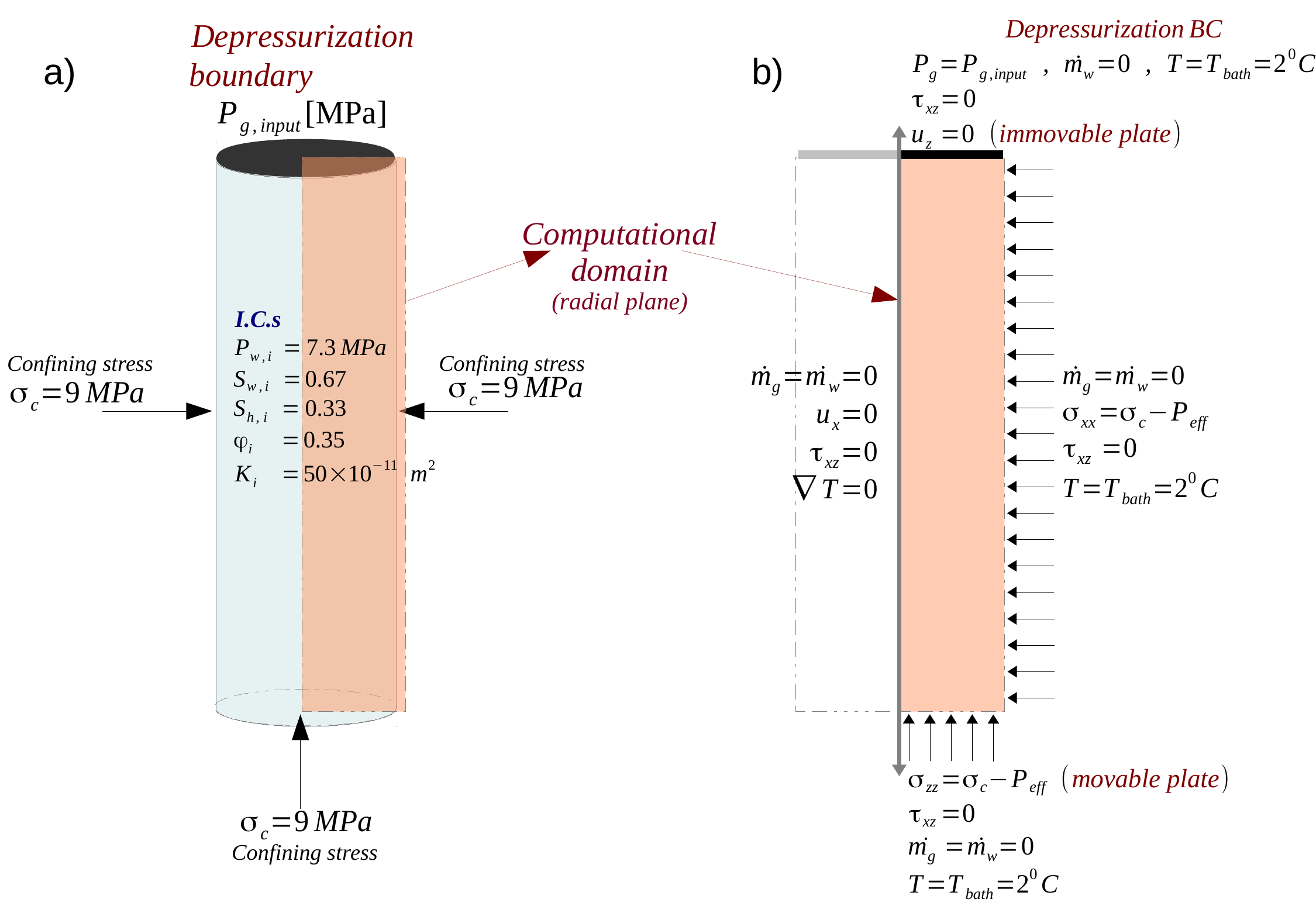}
 \caption{Test setting for the \textbf{depressurization and gas production} period.\newline
	  a) shows the sample and the initial conditions, and b) shows the $2D$ computational domain and the boundary conditions.}
  \label{fig:depressurization_testSetting}
\end{figure*}

\paragraph{b) Depressurization and gas production period}
The schematic of the depressurization test is shown in Fig.\ref{fig:depressurization_testSetting}.
The schematic also shows the initial and boundary conditions.
The simulation is run until $t_{end}=86400$ s (i.e. $1$ day) using a maximum time step size of $120$ s.

% \label{fig:depressurization_testSetting}

\paragraph{Properties and parameters}
The material properties and model parameters chosen for this simulation are listed in Table \ref{table:materialProperties}.
The values of the thermal conductivities, specific heat capacities, dynamic viscosities, and densities for each phase are chosen from standard literature, the references to which are included in the table.
The Brooks-Corey parameters are chosen from the range of typically expected values for sand samples.

The most important properties and parameters relevant to the simulation of the experimental data arise from 
1) the hydrate-phase-change kinetics, and 
2) the poro-elastic behaviour of the hydrate bearing sediments.

\begin{figure}
 \centering
 \includegraphics[scale=0.3]{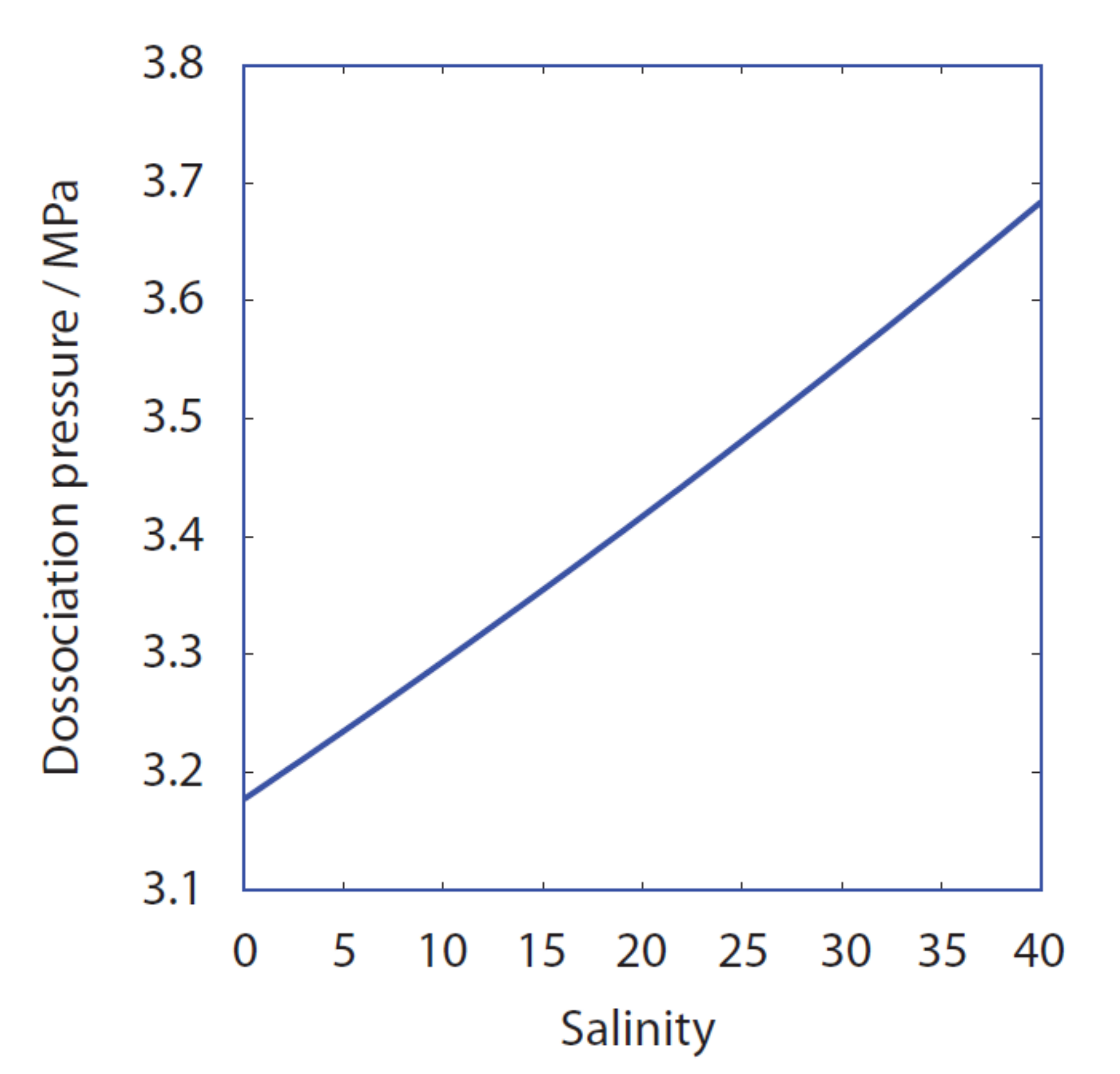}
 \caption{Effect of salinity on hydrate stability curve (at $T_{bath}=2^0C$)}
 \label{fig:hydrateStabilityWithSalinity}
\end{figure}

The hydrate phase-change is modelled by Eqns. (\ref{eqn:GasGenRate},\ref{eqn:WaterGenRate}-\ref{eqn:HeatOfReaction}).
The hydrate-phase equilibrium pressure $P_{e}$ in Eqn. \ref{eqn:GasGenRate} is modelled in accordance with the findings of \cite{Kamath1984}.
For hydrates in pure water, the equilibrium pressure depends only on the temperature.
However, for hydrates in sea water (which is the case for our sample), the equilibrium pressure also depends on the salinity, as shown in Fig.\ref{fig:hydrateStabilityWithSalinity}.
We account for the effect of salinity on the hydrate equilibrium pressure through linear curve fitting on dissociation pressure vs. salinity curve.

% \label{fig:hydrateStabilityWithSalinity}

% \label{fig:formation_results}
  \begin{figure}
    \centering
      \parbox{\textwidth}{a) \textcolor{black}{Average gas pressure $P_g$ in the domain over time.}}
      \vfill
      \includegraphics[scale=0.45]{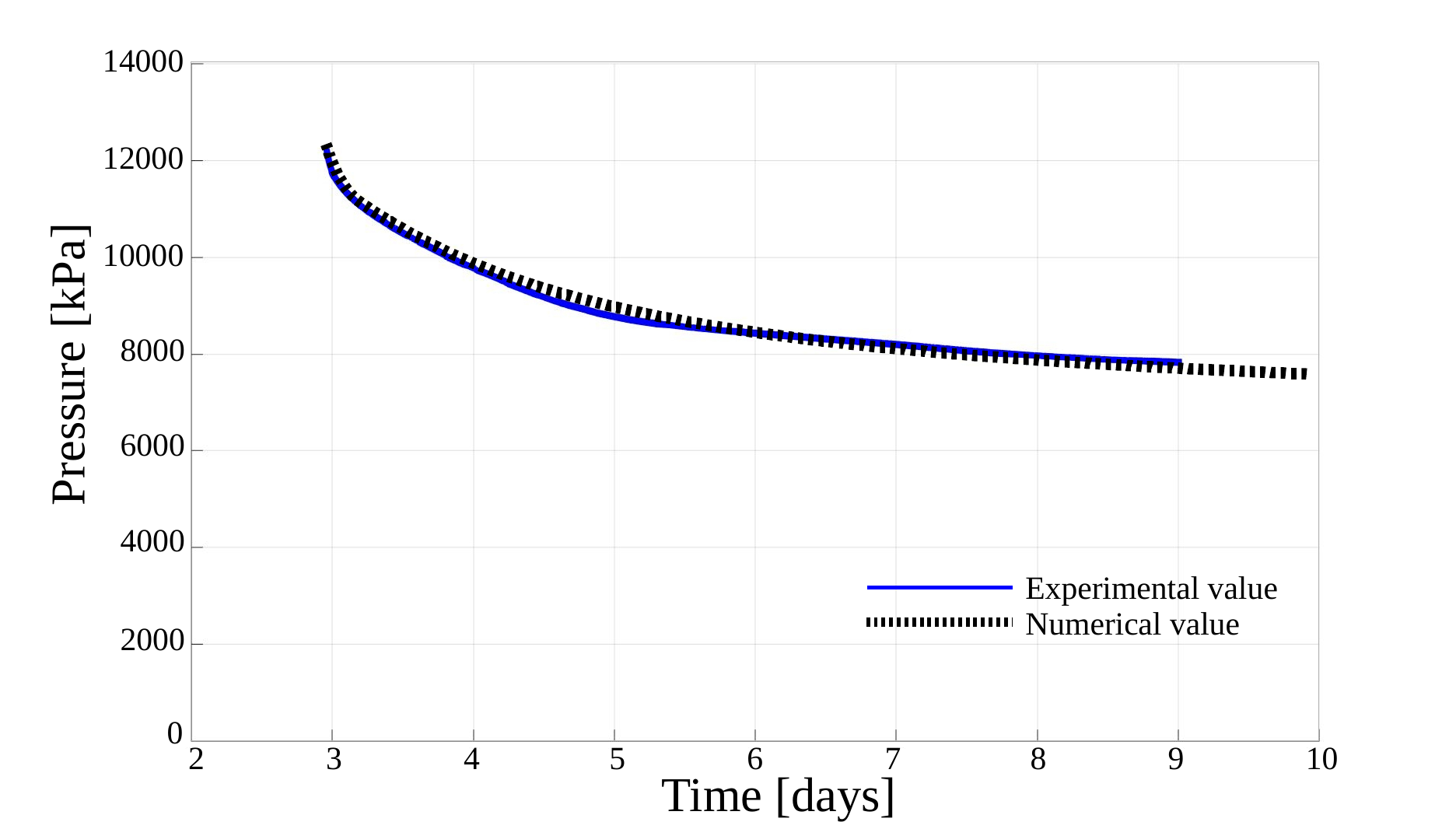}
      \vfill
      \parbox{\textwidth}{b) \textcolor{black}{Average $S_w$ and $S_h$ in the domain over time.}}
      \vfill
      \includegraphics[scale=0.45]{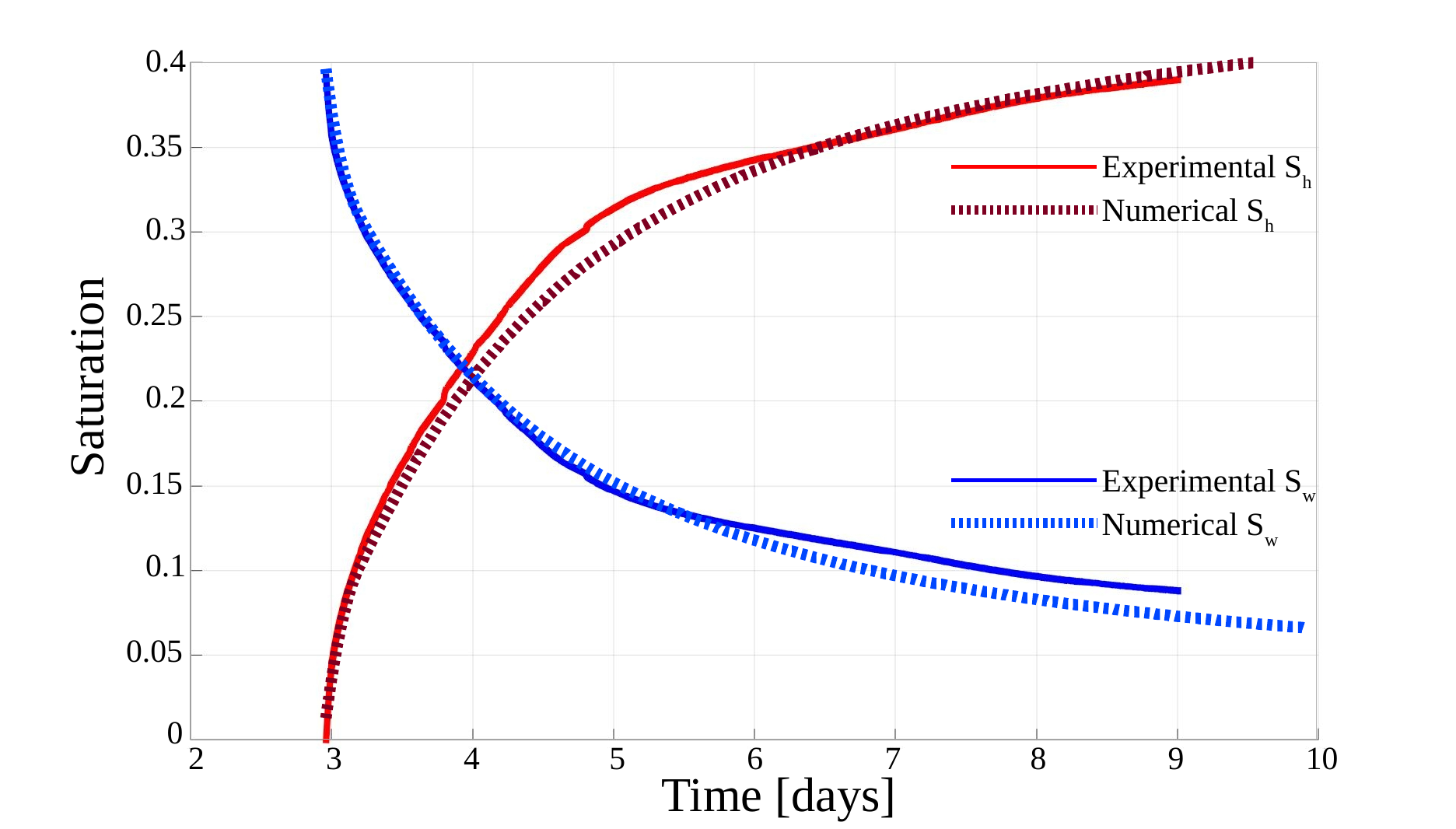}
      \vfill
      \parbox{\textwidth}{c) \textcolor{black}{Total volumetric strain in the sample over time.\newline note: '+' value indicates compression.}}
      \vfill
      \includegraphics[scale=0.45]{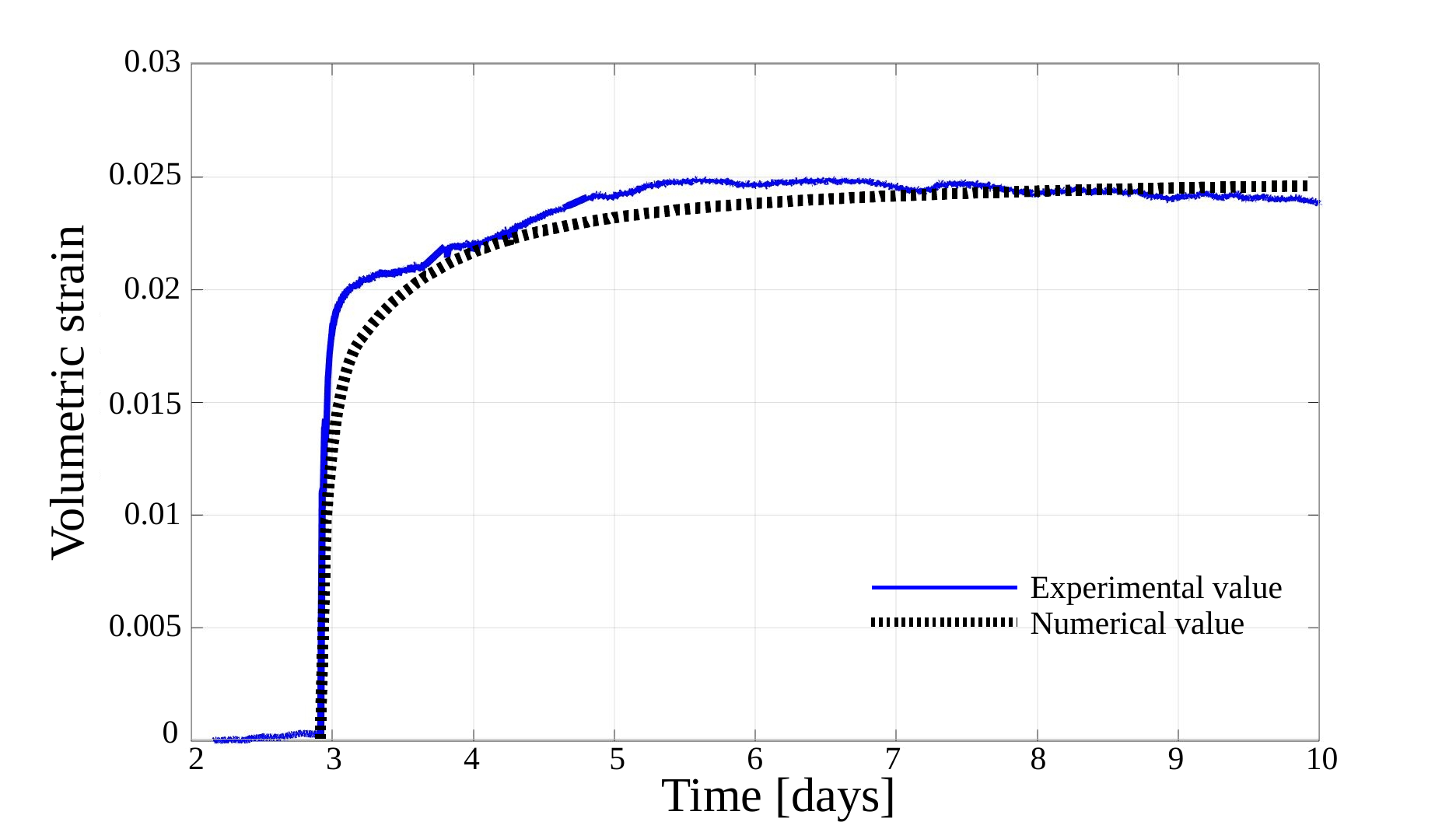}
      \caption{Comparison of the simulation results with the experimental results for the \textbf{gas hydrate formation} period.}
      \label{fig:formation_results}
  \end{figure}

The reaction surface area, $A_{rs}$, in Eqn. \ref{eqn:GasGenRate}, describes the surface area available for the kinetic-reaction, and puts a limit on the mass transfer during hydrate formation and dissociation.
As the hydrate saturation in the pore-space increases, the availability of free surface for hydrate formation to occur decreases, and vice versa. 
Additionally, for hydrate formation, availability of both gas and water in sufficient quantities in the pore-space is a necessary condition. 
This behaviour of $A_{rs}$ is modelled using the parameterization proposed by \cite{SunMohanty2006}.%, given by Eqn. \ref{eqn:Ars}.

The rate of reaction, $k_{reac}$, is a \textit{free parameter} in our simulation which is used to calibrate the hydrate-kinetics model with respect to the experimental data.
In the table we can see that the values of $k_{reac}$, for both hydrate formation as well as dissociation periods, lie well within the range reported in the literature.

The poro-elastic behaviour of the hydrate-bearing sediment is characterized by three parameters, viz., Biot's constant $\alpha_{biot}$, Poisson ratio $\nu_{sh}$, and Young's modulus $E_{sh}$.
Biot's constant is chosen from a range of typically expected values. 
The Poisson's ratio is assumed to be a constant independent of the hydrate saturation following the experimental studies by \cite{Miyazaki2011,Lee2010}.
The Young's modulus is modelled using the parameterization proposed by \cite{SantamarinaRuppel2010}, given by Eqn. \ref{eqn:Esh}.
The Young's modulus $E_{sh}$ is a \textit{free parameter} which are used to calibrate the poro-elasticity model with respect to the experimental data.

\begin{figure}[h]
    \centering
      \parbox{\textwidth}{a) \textcolor{black}{Gas pressure $P_g$ at depressurization boundary (at $z=0$).}}
      \vfill
      \includegraphics[scale=0.45]{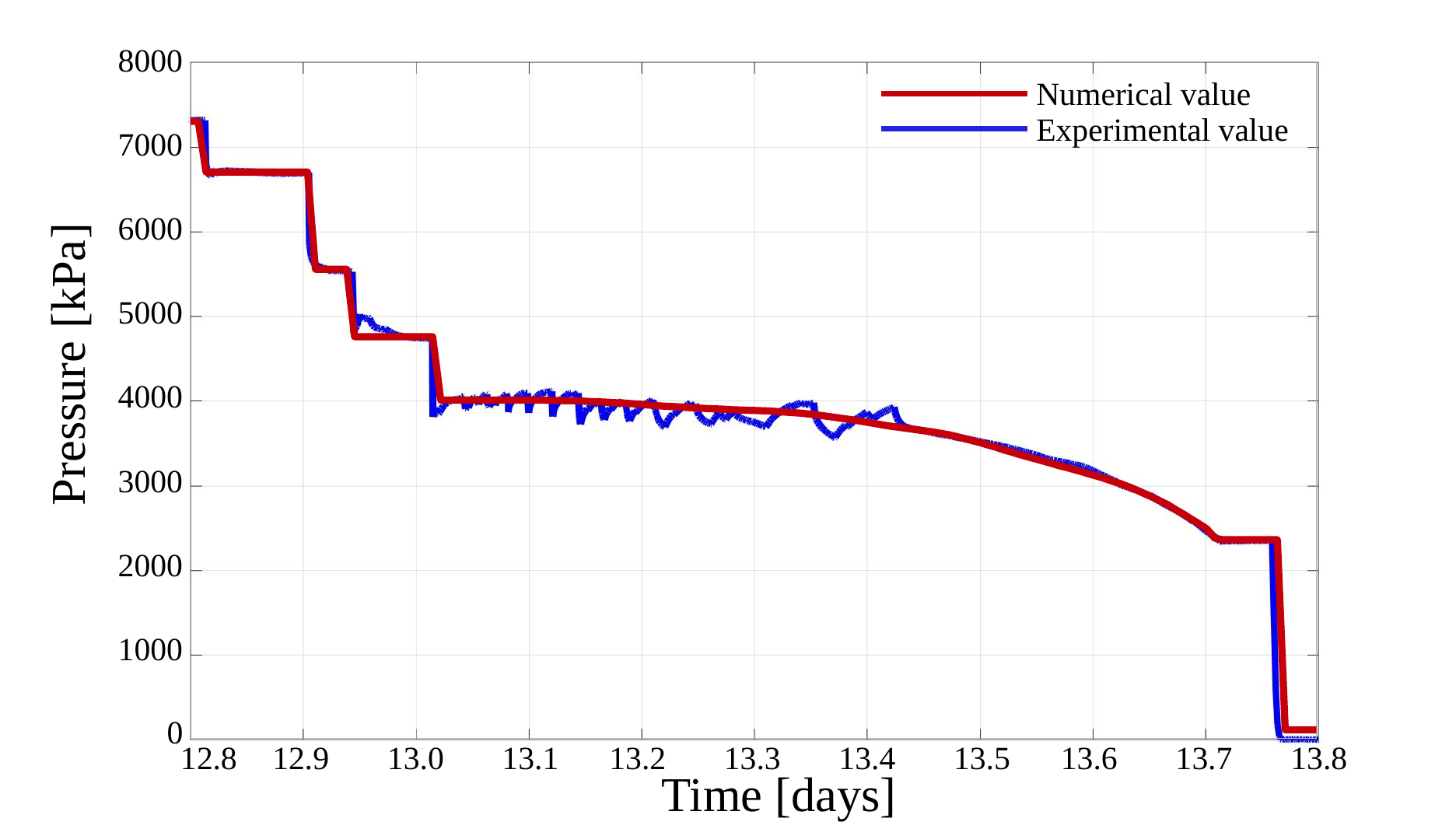}
      \vfill
      \parbox{\textwidth}{b) \textcolor{black}{Cumulative gas production over time.}}
      \vfill
      \includegraphics[scale=0.45]{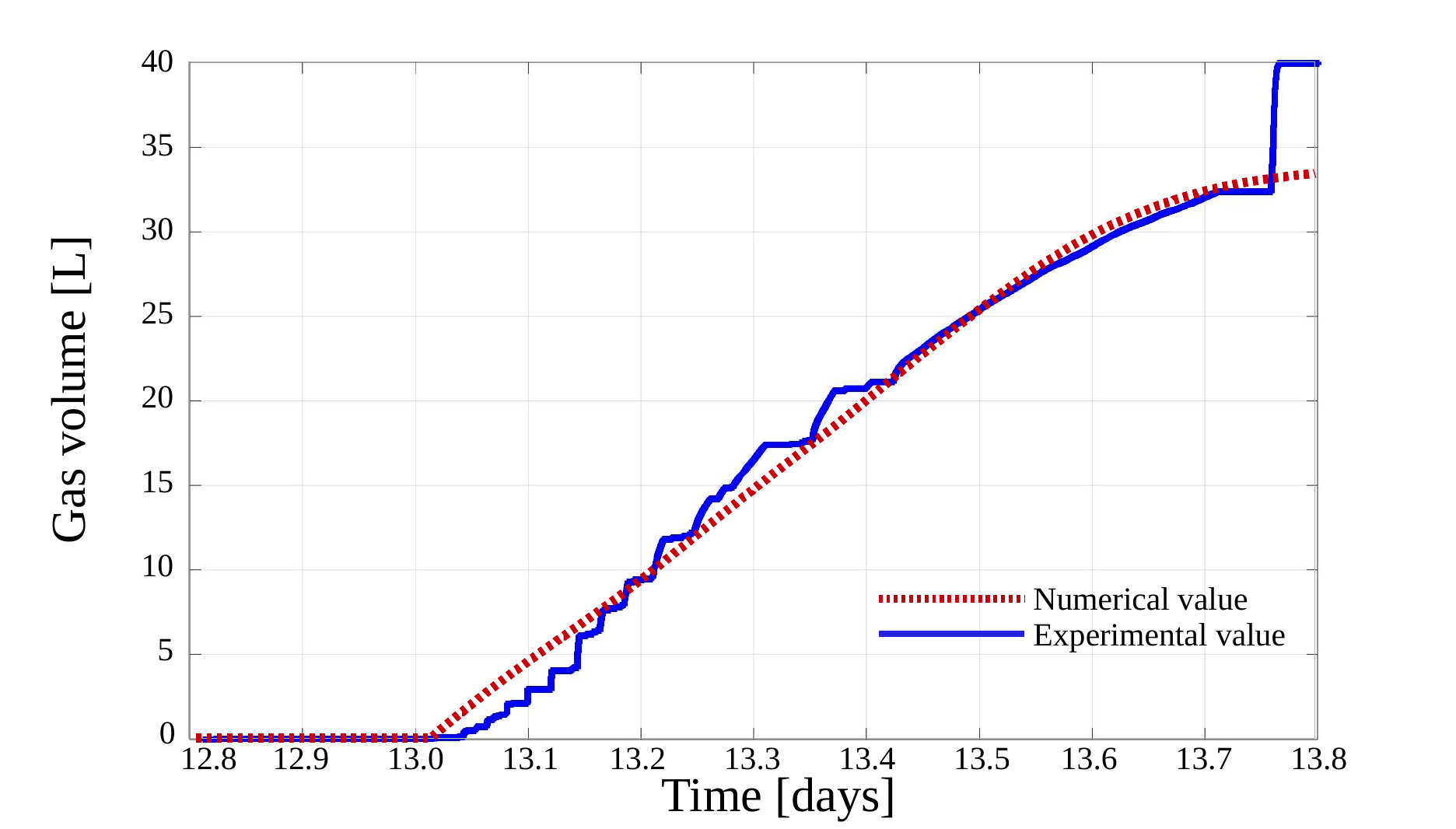}
      \vfill
      \parbox{\textwidth}{c) \textcolor{black}{Total volumetric strain in the sample over time.\newline note: '+' value indicates compression.}}
      \vfill
      \includegraphics[scale=0.45]{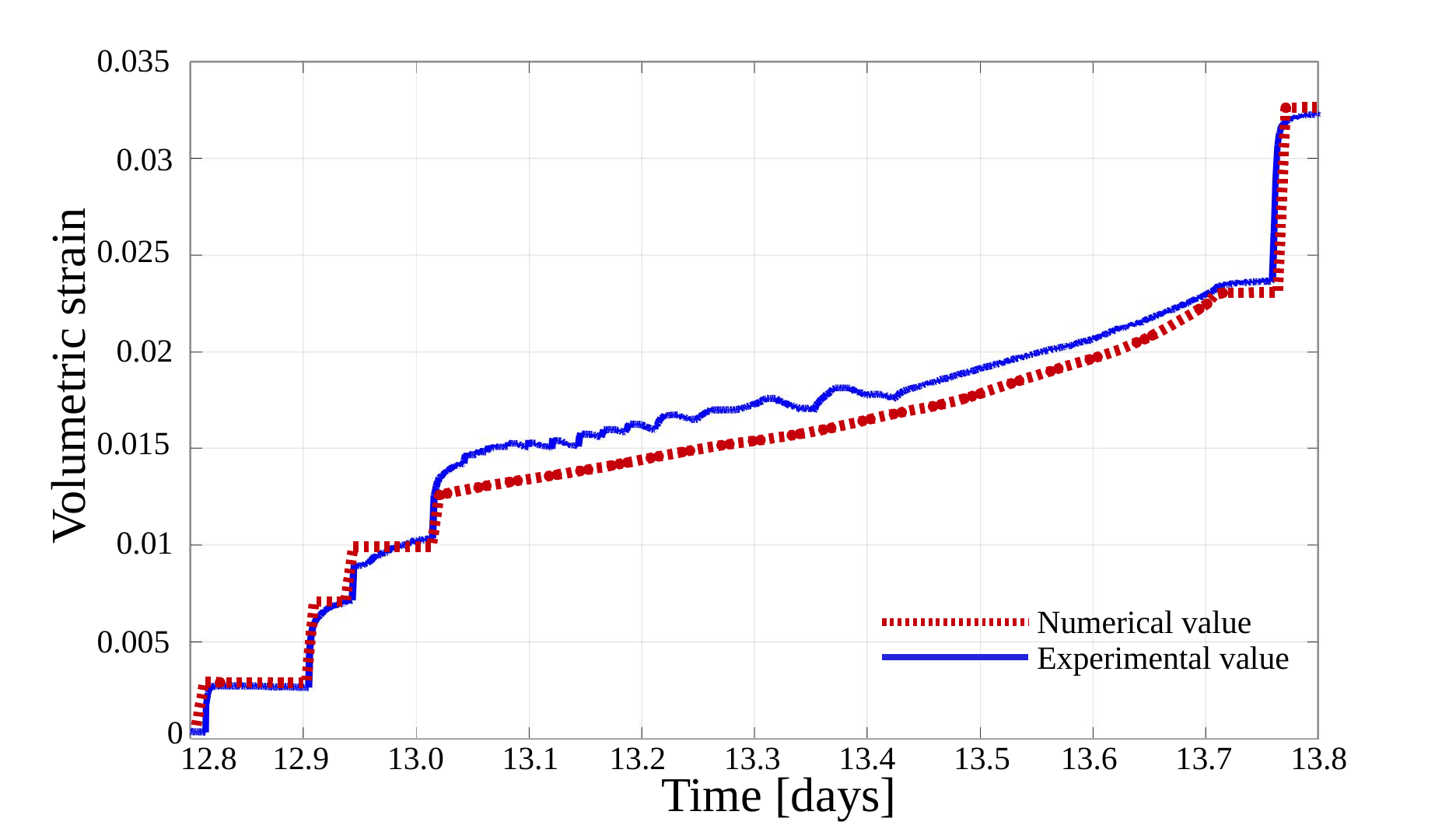}
      \vfill
      \parbox{\textwidth}{d) \textcolor{black}{Temperature in the sample over time.}}
      \vfill
      \includegraphics[scale=0.45]{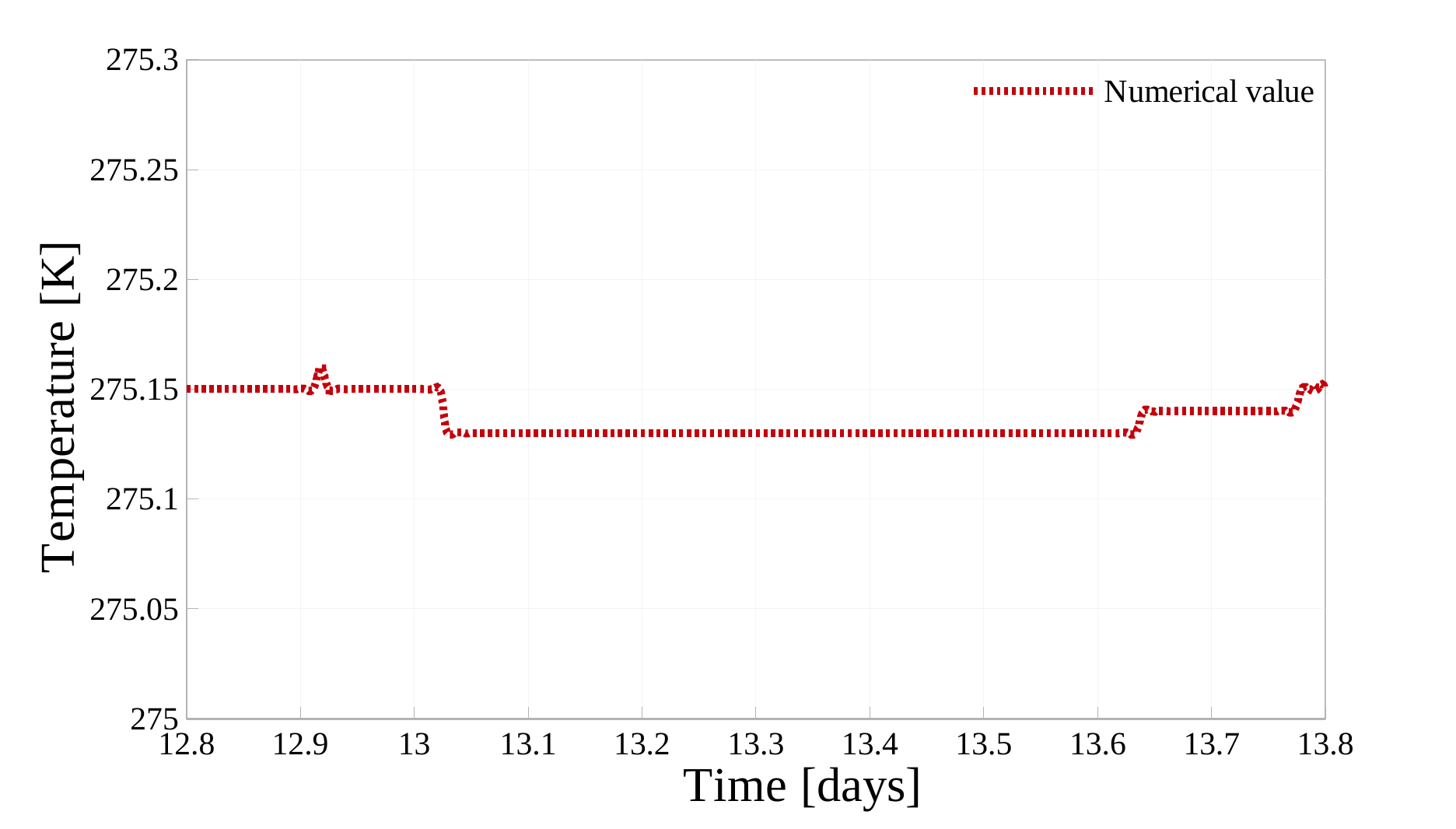}
      \caption{Comparison of the simulation results with the experimental results for the \textbf{depressurization and gas production} period.}
      \label{fig:depressurization_results}
  \end{figure}

\paragraph{Simulation results}
As discussed above, we essentially chose \textit{one} free parameter in kinetics, i.e., $k_{reac}$, and \textit{one} free parameter in linear-elasticity, which is $E_{sh}$,
to calibrate the kinetic and the mechanical models separately. 
With these calibrated models, we simulate numerically the coupled (thermo-chemo)-hydro-geomechanical response of the sand sample in triaxial test-setting using our gas hydrate reservoir model.
The numerical results, together with the corresponding experimental data, are plotted in Fig. \ref{fig:formation_results} for the gas hydrate formation period, 
and in Fig. \ref{fig:depressurization_results} for the hydrate dissociation period.

In the \textbf{gas hydrate formation period}, methane gas in the free pore space is continuously consumed and average bulk gas pressure is decreased. 
(See Fig. \ref{fig:formation_results}-a.)  
Clearly, the rate of gas hydrate formation is not constant. 
In the beginning, after the sample is pressurized at constant isotropic effective stress, gas hydrate formation from free methane gas and pore water is fast, 
but the rate of formation steadily decreases due to mass transfer limitations and shrinking reaction surfaces. 
In accordance to that, after pressurization the gas hydrate saturation increases rapidly and the water saturation decreases proportionally (Fig. \ref{fig:formation_results}-b). 
\textcolor{black}{Note that the reported values of phase saturations are calculated based on initial values and gas pressure measurements.}
\textcolor{black}{The volumetric strain shows a fast positive response during early gas hydrate formation at relatively low gas hydrate saturations, 
and sample stiffness increases at higher gas hydrate saturations (Fig. \ref{fig:formation_results}-c).
The fast volumetric strain response that occurs at constant apparent effective stress results from changes in water saturation and, thus, 
capillary pressure, which is not monitored experimentally, but considered in the numerical simulation.}

During the \textbf{gas hydrate dissociation period}, pressure is decreased step-wise until methane hydrates become unstable at the respective P-/T-conditions. 
Fig. \ref{fig:depressurization_results}-a shows the numerically computed gas phase pressure in the sample.
The gas production is plotted in Fig. \ref{fig:depressurization_results}-b.
With the onset of gas hydrate dissociation after reaching the hydrate stability boundary, 
pressure is maintained at a relatively constant level because hydrate dissociation and gas production equilibrate dependent on experimental and technical conditions.

Volumetric strain during gas hydrate dissociation, plotted in Fig. \ref{fig:depressurization_results}-c, is dependent on effective stress and gas hydrate saturation through the sample stiffness, which decreases with the ongoing gas hydrate dissociation and gas production.

Fig. \ref{fig:depressurization_results}-d shows the numerically computed temperature profile of the sample during dissociation. 
The model predicts that sub-cooling from gas hydrate dissociation is quite small, which is expected since the experiment was performed under isothermal temperature control.

  % \label{fig:depressurization_results}

\section{Discussion and Outlook}
\label{sec:discussionAndOutlook}
In our combined experimental-numerical study, we consider dynamic gas hydrate formation and dissociation in sandy sediment under isotropic compressive loading and show that a simplified coupling concept is capable of reproducing the essential bulk physical behavior, including volumetric strain and gas production. 

The assumption that the soil is the primary load-bearing constituent is \textcolor{black}{central to our model concept}.
This assumption is most reasonable for the pore-filling hydrates with low satuations where the hydrates form by nucleating on sediment grain boundaries and grow freely into pore spaces without entering the pore-throats. 
% In our experiment, we obtain a maximum hydrate saturation of $\approx 0.4$. 
For hydrates formed in partially water saturated sands, as is the case in our experiment, it is well known that the hydrates nucleate preferentially in the pore-throats
and contribute to the sediment stiffness already at low gas hydrate saturations. 
For hydrate saturations between $0.25-0.4$, the hydrates are expected to transition towards a load-bearing habit \cite{Waite2009}.
In our experiment, we obtain a maximum hydrate saturation of $\approx 0.4$.
Our numerical results suggest that in \textit{well consolidated} sands, our assumption of soil forming the primary load-bearing skeleton remains valid even for hydrate saturations which lie in the transition zone between pore-filling and load-bearing habits.
We hypothesize that this is because in \textit{well consolidated} sands, the deformation of hydrates relative to the soil skeleton is quite small in the transition zone.
However, for higher hydrate saturations where hydrates become fully load-bearing, we expect strong limitations to this coupling concept.  
Furthermore, for massive hydrates with saturations exceeding $0.8$, we even expect that the hydrate and the soil phases can no longer be modelled as a single composite phase, and new model concepts are necessary to consistently describe the interface conditions between the hydrate and the soil phase boundaries.

Our experiment was focused on analyzing deformation under variable gas hydrate saturation, and a wide range of effective stress loading, controlled between $1$ to $9$ MPa. 
% To limit the bulk sample deformation and relative grain-to-grain movement, and to avoid critical state conditions and large strain deformation, only isotropic stresses were applied. 
\textcolor{black}{To limit the bulk sample deformation and relative grain-to-grain movement only isotropic stresses were applied.}
Gas hydrates were formed after isotropic consolidation to $1$ MPa using the excess-gas-method (\cite{Chong2016}, \cite{Choi2014}, \cite{Jin2012}, \cite{Priest2009}).
After gas hydrate formation, the remaining gas was fully replaced with seawater. 
After gas-seawater replacement, the sample was equilibrated for approximately $5$ days to allow for gas hydrate alteration before the sample was depressurized. 
A poro-elasticity framework was chosen for describing the mechanical behavior of the sediment in order to minimize the uncertainties arising from unknown mechanical behavior of gas hydrate-bearing sediments. 
The deliberate choice of a simple constitutive law with a limited number of parameters, in contrast to using more complex elasto-plastic modeling approaches, is justified by the design of the experimental test case.
It is important to note that, within the constraints of small-strain deformations and pore-filling hydrates, the coupling concept presented here does not depend on the stress-strain constitutive law as such.
The concept of poro-elasticity is, therefore, sufficient to test the validity of our coupling concept provided that the design of the experiment ensures that the sample deformations remain small and well within elastic limit. 
For large deformations, the coupling concept is not validated so far.
% In our experiment, we obtain a maximum hydrate saturation of $\approx 0.4$ which lies within the range of pore-filling habit.
% For massive gas hydrates with hydrate saturations exceeding $0.4-0.5$, hydrates begin to enter the pore-throats and act as load-bearing solid constituents of the sediment, and it can no longer be assumed that sediment particles alone form the primary load bearnig skeleton. For such hydrates, we expect strong limitations to our coupling concept.

\begin{figure}
  \centering
    \includegraphics[scale=0.45]{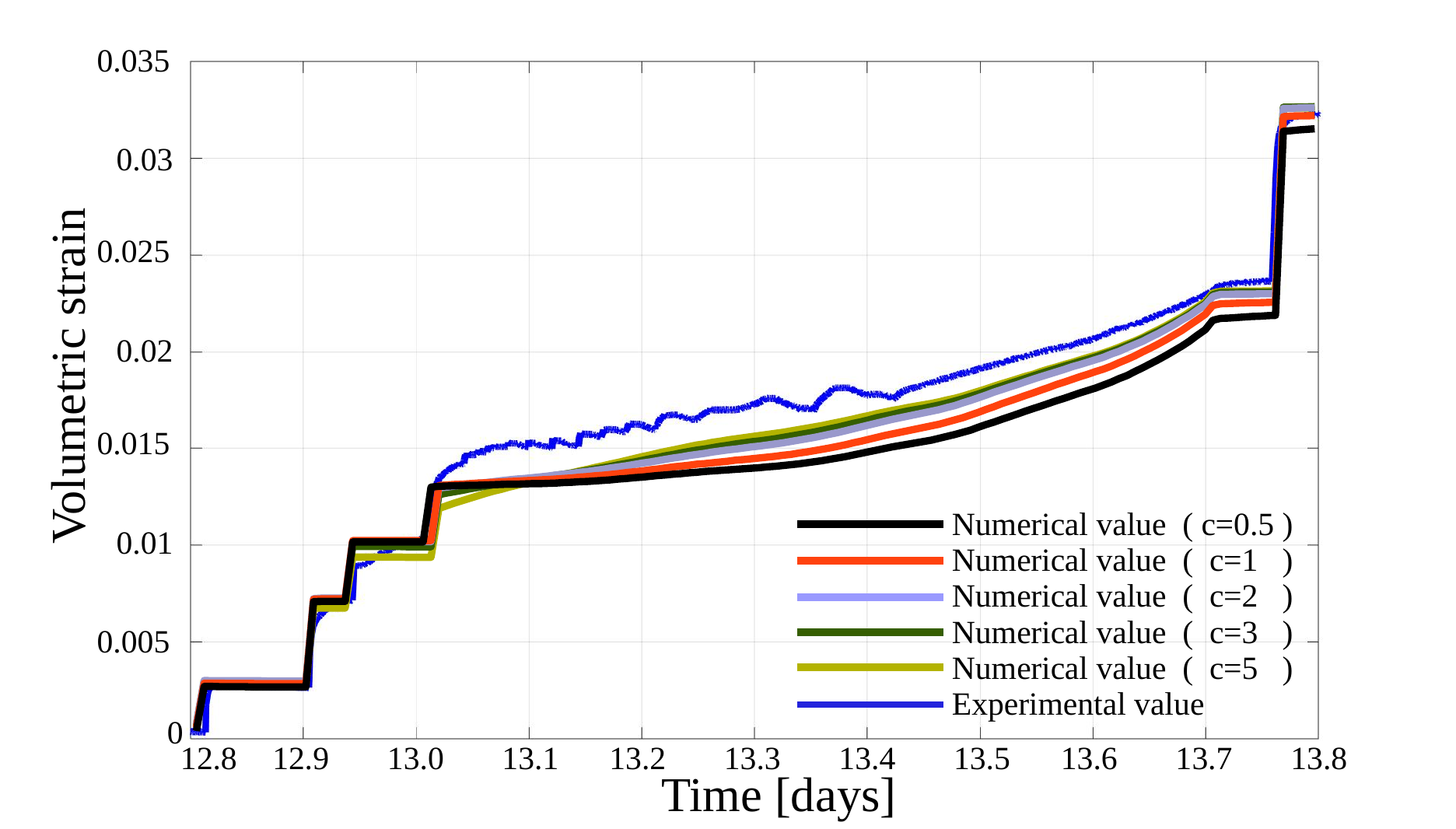}
    \caption{\textcolor{black}{Volumetric strain curves for different functional dependences of $E_{sh}$ on $S_h$ (i.e. $c=1,2,3,5$ ) for the \textbf{depressurization and gas production} period.
    \newline note: '+' value indicates compression.}}
    \label{fig:depressurization_results-volStrain-differentCs}
 \end{figure}

To approximate our experimental data, we treated the kinetic term $k_{reac}$ in the transport block as a fitting parameter.
% , and the fitted parameters are well within the range of values reported in the literature (Lit). 
Similarly, we used fitting of the stiffness model parameters to match the experimental volumetric strain behavior. 
We chose a functional dependence of composite modulus $E_{sh}$ on hydrate saturation $S_h$ as proposed by \cite{SantamarinaRuppel2010} and other authors (\cite{RutqvistMoridis2009}, \cite{KlarSogaNG2010,KlarUchidaSogaYamamoto2013}). 
Knowledge about mechanical stiffness and strength properties of gas hydrate-bearing sediments is still limited. 
Experimental analysis of mechanical properties is problematic, because it is equally important to control effective stress conditions and phase saturations, and there are no test procedures available to guarantee homogeneous gas hydrate saturations and full water saturation. 
Further, mechanical properties are strongly dependent on gas hydrate-sediment fabrics and formation procedures, and effects from dynamic changes in gas hydrate saturation, distributions and alterations of gas hydrate-sediment fabrics needs further investigation and development of novel test procedures (\cite{DeusnerICEGT2016}) particularly for dynamic test scenarios. 
Overall, the calibrated values for $E_{sh}$ from our study are in accordance with earlier experimental and numerical studies, which reported Young’s modulus or secant stiffness in a
wide range of approximately $100$ to $400$ MPa for relevant gas hydrate concentrations (\cite{Brugada2010}, \cite{Miyazaki2010}, \cite{LeeSantamarinaRuppel2010}, \cite{Yun2007}). 
Although the composite modulus $E_{sh}$ was treated as a free fitting parameter and initialized based on apparent stress-strain behavior during the intervals of known and constant gas hydrate saturations, physically meaningful values for individual modulii $E_s$ and $E_h$ were obtained. 
$E_s$ for the sediment without gas hydrate reflected stiffness behavior typical for loose soil during gas hydrate formation while the sample was normally consolidated at low effective stress. 
\textcolor{black}{The results from the numerical simulation suggest that an apparent step-like increase in bulk composite
stiffness (i.e., the change of apparent $E_{sh}$ from $132$ MPa to $183$ MPa) occurred during the time interval
between the completion of gas hydrate formation and the start of depressurization. 
We assume that this response was caused by the high transient effective stress and the composite sediment consolidation, 
which could not be avoided during the gas-water exchange. 
In order to accomplish the replacement of gas with water sufficiently fast, and to minimize gas hydrate dissociation during the short time interval of gas-water exchange, 
the confining stress instead of the apparent effective stress was controlled at a constant value during gas-water exchange. 
In fitting the model to the experimental data, we adjusted both $E_s$ and $E_h$ rather than constraining the effect of consolidation to one of the parameters a priori. 
However, the validity of this assumption needs to be further investigated using advanced geotechnical and microstructural analyses. 
Furthermore, it needs to be considered that the poro-elasticity model adopted for testing our coupled numerical simulation scheme does not explicitly consider effective stress-dependent changes in the modulii $E_s$ and $E_h$,
which could also contribute to the apparent differences in $E_{sh}$ during gas hydrate formation and dissociation periods.}
\textcolor{black}{In our simulation,} the composite modulus $E_{sh}$ depends almost linearly on $S_h$ during gas hydrate formation, while during the hydrate dissociation period the dependence of $E_{sh}$ on $S_h$ is smaller. 
\textcolor{black}{Fig. \ref{fig:depressurization_results-volStrain-differentCs}} shows the volumetric-strain plotted over time for the depressurization period for different functional dependences of $E_{sh}$ on $S_h$ (i.e. $c = 0.5, 1, 2, 3, 5$ ). 
Our simulation results indicate that an exponent $c = 3$ is a reasonable approximation. 
\cite{SantamarinaRuppel2010} suggest that $S_h$ tends to be raised to a power larger than $1$, which reduces the impact on stiffness at low gas hydrate saturations relative to that for high gas hydrate saturations. 
Since gas hydrates formed using the excess-gas-method are predominantly located in the pore throats rather than in the free pore space, the linear and relatively stronger dependence of $E_{sh}$ on $S_h$ during formation appears reasonable. 
The weaker dependence of $E_{sh}$ on gas hydrate saturation during dissociation is also reasonable, since after exchanging gas with water in the pore space, gas hydrate-sediment fabrics were allowed to alter, and also during dissociation the grain-scale hydrate-sand structure is necessarily changed. 
Thus, our results clearly show that dynamic structural transitions in gas hydrate-bearing sediments during gas hydrate formation, aging and dissociation can have substantial effects on sediment mechanicalproperties. 
Further combined experimental-numerical studies with the objective to simulate the geomechanical effects from such dynamic changes in gas hydrate-sediment fabrics are currently
ongoing.

\newpage
% mathematical model table
\begin{longtable}{L{3.5cm} L{10cm} C{1cm}}
\caption{Summary of the mathematical model}
\label{table:modelSummary}\\
\specialrule{.1em}{.05em}{.05em}
\endfirsthead
\multicolumn{2}{c}%
{\tablename\ \thetable\ -- Summary of the mathematical model (\textit{continued})} & \scalebox{0.8}{eqn. no.}\\
\hline
\endhead
\hline \multicolumn{3}{r}{\textit{Continued on next page}} \\
\endfoot
\specialrule{.1em}{.05em}{.05em}
\endlastfoot
% \hline
  \multicolumn{2}{c}{\textbf{Governing equations}} & \scalebox{0.8}{eqn. no.} \\ \hline
  Mass balance for each mobile component $\kappa=CH_4,H_2O$	& $\sum\limits_{\alpha} \left[ \partial_t \left( \phi \ \rho_{\alpha} \  \chi_{\alpha}^{\kappa} \ S_{\alpha} \right) \right]
							  \ + \ \sum\limits_{\alpha} \left[ \nabla\cdot\left( \phi \ \rho_{\alpha} \  \chi_{\alpha}^{\kappa} \ S_{\alpha} \ {\mathbf v_{\alpha,t}} \right) \right]
							    = \sum\limits_{\alpha} \left[ \nabla\cdot\left(\phi \ S_{\alpha} \ {\bf J}_{\alpha}^{\kappa}\right) \right]  
							  \ + \ \dot g^{\kappa}
							  \ + \ \sum\limits_{\alpha} \dot q_{m_{\alpha}}^{\kappa}
							  \ + \ \dot S_{ext}^{\kappa}$ 					& (\Eqn{eqn:MassBal_g}),(\Eqn{eqn:MassBal_w})\\ 
% \\
  Mass balance for hydrate phase			& $\partial_t \left(\phi \ \rho_{h} \  S_{h} \right)
							  \ + \ \nabla\cdot\left(\phi \ \rho_{h} \  S_{h} \ {\mathbf v_{h,t}} \right)
							  \ = \ \dot g^{h}$ 						& (\Eqn{eqn:MassBal_h})\\
% \\
  Mass balance for sand phase				& $\partial_t \left[(1-\phi) \ \rho_{s} \right]
							  \ + \ \nabla\cdot\left((1-\phi) \ \rho_{s} \ {\mathbf v_{s}} \right)
							  \ = \ 0$ 							& (\Eqn{eqn:MassBal_s})\\ 
% \\
  Momentum balance for mobile phases $\alpha=g,w$	& ${\bf v}_{\alpha}=-K \ \frac{k_{r\alpha}}{\mu_{\alpha}}\left(\nabla P_\alpha - \rho_{\alpha} \ {\bf g}\right)$ \quad (Darcy's Law) 
															& (\Eqn{eqn:Darcy_g}),(\Eqn{eqn:Darcy_w}) \\
% \\
  Momentum balance for 				 	& $\nabla \cdot \tilde{\sigma} + \rho_{m} \ {\bf g} = 0$ 	& (\Eqn{eqn:MomentumBal_sh})\\
  composite - solid					& where, $\rho_{m}$ is the bulk density given by 
							  $\ \rho_{m} = \sum\limits_{\beta}\left(\phi \ S_{\beta} \ \rho_{\beta}\right) + \left(1-\phi\right) \rho_{s}$ 
															& \\
% \\
  Energy balance					& $\partial_t \left[ \left( 1-\phi \right) \rho_s u_s + \sum_{\beta} \left( \phi \ S_{\beta} \ \rho_{\beta} \ u_{\beta} \right)  \right]
							   \ + \ \sum_{\alpha} \left[ \nabla \cdot \left( \phi \ \rho_{\alpha} \ \chi_{\alpha}^{\kappa} \ S_{\alpha} \ {\bf v }_{\alpha,t} \ h_{\alpha} \right) \right]
							   \ = \ \nabla \cdot k^c_{eff} \nabla T 
							   \ + \ \dot Q_h 
							   \ + \ \sum\limits_{\alpha} \left( \dot q_{m_{\alpha}}^{\kappa} \  h_{\alpha} \right)$ 	
															& (\Eqn{eqn:EnergyBal})\\
							& where, 							&\\
							& $  k^c_{eff} = \left( 1 - \phi \right) k^c_s + \sum_{\alpha}\sum_{\kappa} \left( \phi \ \chi_{\alpha}^{\kappa} \ S_{\alpha} \ k_{\alpha}^c \right) 
							     + \phi \ S_{h} \ k_h^c $ 					& \\
							& $  h_{\alpha} = \int_{T_{ref}}^T Cp_{\alpha}\ dT $ 		&\\
							& $  u_{\gamma} = \int_{T_{ref}}^T Cv_{\gamma}\ dT $ 		&\\[0.75em]
\hline
  \multicolumn{3}{c}{\textbf{Closure relationships}} \\ \hline
  Relationship between phase pressures 	& $P_g - P_w = P_c\left( S_{we}\right)$ 					&  (\Eqn{eqn:Closure1})	\\
  Summation relationships 		& $\sum\limits_{\beta} S_{\beta}=1$  						& (\Eqn{eqn:Closure2})	\\
					& $\forall \alpha$ : $\sum\limits_{\kappa}\chi_{\alpha}^{\kappa}=1$		& (\Eqn{eqn:Closure3}),(\Eqn{eqn:Closure3}) \\
\hline
  \multicolumn{3}{c}{\textbf{Constitutive relationships}} \\ \hline
  \multicolumn{3}{l}{1. \underline{Vapor-liquid equilibrium}} \\
  \multicolumn{2}{p{13.5cm}}{Using Henry's Law and Raoult's Law for ideal gas - liquid solutions,}
																&\\
  \multicolumn{2}{l}{For dissolved methane:	\quad	$\ \chi_w^{CH_4} = H (T) \ \chi_g^{CH_4} \ P_g$} 			&\\
  \multicolumn{2}{l}{For water vapor:		\qquad\qquad  $\chi_g^{H_2 O} = \chi_w^{H_2 O} \ \frac{P^{sat}_{H_2 O} (T) }{P_g} $} 
																&\\
  \multicolumn{2}{p{13.5cm}}{where, $H(T)$ is the Henry's constant for methane dissolved in water calculated using the empirical relation from NIST database (\cite{NISTChemHenryConst})
, and $P^{sat}_{H_2 O}$ is saturated water vapour pressure calculated using Antoine's equation. }
																&\\ 
  \multicolumn{3}{p{14.5cm}}{2. \underline{Diffusive mass-transfer flux}}\\[0.25em]
  Fick's law:	&${\bf J}_{\alpha}^{\kappa} = - \tau D^{\alpha} \left( \rho_{\alpha} \nabla \chi_{\alpha}^{\kappa} \right)$ 	&\\
		&where, $D^{\alpha}$ is the binary diffusion coefficient.							 
		$D^g$ is estimated using the empirical relationship proposed by \cite{StatteryBird1958} for low density binary $CH_4-H_2O$ system.
		$D^w$ is estimated using the Wilke-Chang correlation (\cite{Himmelblau1964}) for dilute associated liquid mixtures is used. &\\
  \multicolumn{3}{p{14.5cm}}{3. \underline{Hydrate phase change kinetics}}																\\
  \multicolumn{2}{p{13.5cm}}{Non-equilibrium phase change of methane hydrate is modeled by the Kim-Bishnoi kinetic model (\cite{KimBishnoi1987}).} 		&					\\
  \multicolumn{2}{p{13.5cm}}{Gas generation rate		\quad\quad\qquad $ \dot g^{CH_4} = k_{reac}\ M_{CH4}\ A_{rs}\ \left( P_e - f_g \right) $} 	& (\Eqn{eqn:GasGenRate})		\\
  & where, 
    $k_{reac}$ is the rate of kinetic phase-change,
    $P_e$ is the hydrate equilibrium pressure, modelled as (\cite{Kamath1984}),											&					\\
  & \qquad $P_e = A_1 \ \exp\left( A_2 - \frac{A_3}{T} \right)$ , 												& (\Eqn{eqn:Peq})			\\
  & $f_g$ is the gas fugacity calculated using the Peng-Robinson EoS for methane gas, and,
    $A_{rs}$ is the specific reaction surface area modelled as $A_{rs} = \Gamma_{r} A_{s}$, where, $A_s$ is the total surface area and $\Gamma_r$ is the ratio of the active reaction surface to the total surface area.
    $A_s$ and $\Gamma_r$ are modelled using the correlations proposed by \cite{Yousif1991} and \cite{SunMohanty2006}, respectively.									\\
%     leading to, 																		&					\\
%   & \qquad $A_{rs} = \left\lbrace \genfrac{}{}{0pt}{}{\text{\scalebox{1.25}{$\phi S_h \sqrt{\frac{\phi_{eff}^3}{2K}} \text{ \qquad\qquad, if } (P_e - f_g ) \geq 0$}} }
% 					      {\text{\scalebox{1.25}{$\left( S_g  S_w  S_h \right)^{\frac{2}{3}} \sqrt{\frac{\phi_{eff}^3}{2K}} \text{\quad, if } (P_e - f_g ) < 0$}} }
% 	      \right.$ 																		& (\Eqn{eqn:Ars})			\\
  \multicolumn{2}{p{13.5cm}}{Water generation rate		\quad\qquad $ \dot g^{H_2O} = \dot g^{CH_4} \ N_{Hyd} \ \frac{M_{H_2O}}{M_{CH_4}} $} 		& (\Eqn{eqn:WaterGenRate})		\\
  \multicolumn{2}{p{13.5cm}}{Hydrate consumption rate	\quad $ -\dot g^{Hyd}  = \dot g^{CH_4}\ \frac{M_{Hyd}}{M_{CH_4}} $} 					& (\Eqn{eqn:HydGenRate})		\\
  \multicolumn{2}{p{13.5cm}}{Heat of hydrate dissociation	\ $ \dot Q_h = \frac{-\dot g^{Hyd}}{M_{Hyd}}\left( B_1 - \frac{B_2}{T} \right) $} 		& (\Eqn{eqn:HeatOfReaction})		\\ 
  \multicolumn{3}{p{14.5cm}}{4. \underline{Properties of the fluid-matrix interaction}}															\\
  Capillary pressure		& $P_c = P_{c0}\cdot f^{Pc}_{S_h}\left(S_h\right) \cdot f^{Pc}_{\phi}\left(\phi\right)$ 					& (\Eqn{eqn:Pc})			\\
				& where, 
				  $P_{c0}$ is capillary pressure for un-deformed, un-hydrated solid matrix, given by Brooks-Corey relationship, 		& 					\\
				& \qquad $P_{c0} = P_{entry} \ S_{we}^{\ -1 / \lambda_{BC}}$ , 									& (\Eqn{eqn:Pc_BC})			\\
				& and, $f^{Pc}_{S_h}\left(S_h\right)$ and $f^{Pc}_{\phi}\left(\phi\right)$ are scaling factors to account for the effects of 
				$S_h$ (\cite{Rockhold2001}) and $\phi$ (Civan's power-law correlation \cite{Civan2000}) respectively,				&					\\
				& $f^{Pc}_{S_h} = \left(1-S_h\right)^{-\frac{3 \lambda_{BC} - 1}{3 \lambda_{BC}}}$  
				  \ and, $\ f^{Pc}_{\phi} = \frac{\phi_0}{\phi} \left( \frac{1-\phi}{1-\phi_0}\right)^{2}$ 					& (\Eqn{eqn:SF_Pc})			\\
  Intrinsic permeability	& $K = K_{0}\cdot f^{K}_{S_h}\left(S_h\right) \cdot f^{K}_{\phi}\left(\phi\right)$ 						& (\Eqn{eqn:K})				\\
				& where, 
				  $K_{0}$ is intrinsic permeability of the un-deformed, un-hydrated solid matrix,
				  and, $f^{K}_{S_h}\left(S_h\right)$ and $f^{K}_{\phi}\left(\phi\right)$ are scaling factors to account for the effects of 
				  $S_h$ (\cite{Rockhold2001}) and $\phi$ (\cite{Civan2000}) respectively, 							&					\\
				& \qquad $f^{K}_{S_h} = \left(1-S_h\right)^{\frac{19}{6}}$ \ and, $\ f^{K}_{\phi} = \frac{\phi}{\phi_0} \left( f^{Pc}_{\phi} \right)^{-2}$ 	& (\Eqn{eqn:SF_K})	\\
  Relative permeabilities	& Relative permeabilities of the mobile phases are modelled using Brooks-Corey model in conjunction with the Burdine theorem (\cite{Burdine1953}), as 			&\\
				& $k_{rw} = S_{we}^{\frac{2+3\lambda_{BC}}{\lambda_{BC}}} $ \ and,  
				  $k_{rg} = \left( 1-S_{we} \right)^{2} \left( 1 - S_{we}^{\frac{2+\lambda_{BC}}{\lambda_{BC}}}\right)$ where, $S_{we}=\frac{S_w}{1-S_h}$ 	& (\Eqn{eqn:relPerms})	\\
  Hydraulic tortuosity		& $\tau = \phi^n \quad$ where, $1 \leq n \leq 3$ 												& (\Eqn{eqn:tortuosity})\\ 
  \multicolumn{3}{p{14.5cm}}{5. \underline{Poro-elasticity}}\\
  Effective stress principle		& Effective stress concept introduced by \cite{Terzaghi1925} and modified by \cite{Biot1941} is used.	&\\
					& $\tilde{\sigma} = \tilde{\sigma}^{\prime} + \alpha_{biot} \left(\frac{S_g P_g + S_w P_w }{S_g + S_w} \right) \tilde{I}$ 		& (\Eqn{eqn:effectiveStress})	\\
					& where, $\tilde{\sigma}$ is the total stress acting on the bulk porous medium, and
					  $\tilde{\sigma}^{\prime}$ is the effective stress acting on the composite skeleton. 							&				\\
					& $\alpha_{biot}$ is the Biot-Willis constant (\cite{BiotWillis1957}).  								&				\\
  Linear elastic law			& $\tilde{\sigma}^{\prime} = 2 \ G_{sh} \ \tilde{\epsilon} + \lambda_{sh} (tr \ \tilde{\epsilon})\ \tilde{I}$ 				& (\Eqn{eqn:stress-strain})	\\
					& where, $\tilde{\epsilon}$ is the linearized strain, given by $\quad \tilde{\epsilon} = \frac{1}{2}\left( \nabla {\bf u} + \nabla^T {\bf u}\right)$ 	
																						& (\Eqn{eqn:strain})		\\
					& and, $G_{sh}$ and $\lambda_{sh}$ are the Lame's parameters. 										&				\\
% 					  In terms of Young's modulus $E_{sh}$, and Poisson's ratio $\nu_{sh}$, 								&				\\
% 					& $G_{sh} = \frac{E_{sh}}{2\left(1 + \nu_{sh}\right)} \ , \ \lambda_{sh} = \frac{E_{sh}\ \nu_{sh}}{\left(1 + \nu_{sh}\right)\left(1 - 2\ \nu_{sh}\right)}$ 
% 																						& (\Eqn{eqn:LameParams}) 	\\
  Young's modulus			& Young’s modulus $E_{sh}$ is modelled using the parameterization proposed by \cite{SantamarinaRuppel2010}:
					$E_{sh} = E_{s} + S_{h}^c \ E_{h}$ 													& 				\\
					& where, $E_s$ and $E_h$ are the Young's modulus of hydrate-free sand and hydrate, respectively.					&(\Eqn{eqn:Esh})		\\
\hline
\end{longtable}
\addtocounter{equation}{\theEqnno}

% \newpage 
% material properties table
% \onecolumn
\begin{longtable}{L{0.35cm} R{4.cm} L{4.cm} L{2.5cm} C{3cm}}
\caption{Material properties and model parameters}
\label{table:materialProperties}\\
\specialrule{.1em}{.05em}{.05em}
\endfirsthead
\multicolumn{4}{c}%
{\tablename\ \thetable\ -- Material properties and model parameters (\textit{continued})} \\
\hline
\endhead
\hline \multicolumn{5}{r}{\textit{Continued on next page}} \\
\endfoot
\specialrule{.1em}{.05em}{.05em}
\endlastfoot
      \multicolumn{4}{c}{\textbf{Thermal conductivities}} 		& Ref.\\
      $k^c_g$				& \multicolumn{2}{R{8.cm}}{$	- 0.886\times 10^{-2} 
									+ 0.242\times 10^{-3} T 
									- 0.699\times 10^{-6} T^2 
									+ 0.122\times 10^{-8} T^3$}	& $W\cdot m^{-1}\cdot K^{-1}$	& \cite{Roder1985}\\
      $k^c_w$				& \multicolumn{2}{c}{$0.3834\ ln(T) - 1.581$}			& $W\cdot m^{-1}\cdot K^{-1}$	& \cite{IAPWS1997}\\
      $k^c_h$				& \multicolumn{2}{c}{$2.1$}					& $W\cdot m^{-1}\cdot K^{-1}$	& \cite{Sloan2007Book}\\
      $k^c_s$				& \multicolumn{2}{c}{$1.9$}					& $W\cdot m^{-1}\cdot K^{-1}$	& \cite{Esmaeilzadeh2008}\\
% \\
%       
      \multicolumn{4}{c}{\textbf{Specific heat capacities}}		\\ 
      $Cp_g$				& \multicolumn{2}{r}{$\Delta Cp^{res}_{g}\left( 1238 + 3.13 T + 7.9\times10^{-4} T^2\right.$} &\\
					& \multicolumn{2}{r}{$\left. - 6.86\times10^{-7} T^3\right)$}
										& $J\cdot kg^{-1}\cdot K^{-1}$ & \cite{PengRobinson1976,Esmaeilzadeh2008}\\
      $Cv_g$				& \multicolumn{2}{c}{$Cp_g+R_{CH_4}$}	& $J\cdot kg^{-1}\cdot K^{-1}$ 	& \\
      $Cp_w$				& \multicolumn{2}{c}{$4186$}		& $J\cdot kg^{-1}\cdot K^{-1}$ 	& \cite{IAPWS1997}\\
      $Cv_w$				& \multicolumn{2}{c}{$Cp_w+R_{H_2O}$}	& $J\cdot kg^{-1}\cdot K^{-1}$ 	& \\
      $Cv_h$				& \multicolumn{2}{c}{$2700$}		& $J\cdot kg^{-1}\cdot K^{-1}$ 	& \cite{Sloan2007Book}\\
      $Cv_s$				& \multicolumn{2}{c}{$800$}		& $J\cdot kg^{-1}\cdot K^{-1}$ 	& \cite{Esmaeilzadeh2008}\\ 
% \\
%       
      \multicolumn{4}{c}{\textbf{Dynamic viscosities}} \\
      $\mu_g$	& \multicolumn{2}{c}{$10.4 \times 10^{-6} \left(\frac{273.15 + 162 }{T+162}\right) \left(\frac{T}{273.15}\right)^{1.5}$}	& $Pa\cdot s$	& \cite{MethanePropertiesFriend1989}\\
      $\mu_w$	& \multicolumn{2}{r}{$0.001792 \ \exp\left[ - 1.94 - 4.80 \left(\frac{273.15}{T}\right)\right.$}	\\
		& \multicolumn{2}{r}{$\left. + 6.74 \left(\frac{273.15}{T}\right)^2\right]$}								& $Pa\cdot s$	& \cite{IAPWS1997}\\
% \\
% 
      \multicolumn{4}{c}{\textbf{Densities}} \\
      $\rho_g$	& \multicolumn{2}{c}{$\frac{P_g}{z R_g T}$}		& $kg\cdot m^{-3}$ 	& \cite{PengRobinson1976}\\
      $\rho_w$	& vapour:		&$0.0022\frac{P_g}{T}$		& $kg\cdot m^{-3}$ 	& \cite{IAPWS1997}\\
		& liquid:		&$1000$				& $kg\cdot m^{-3}$ 	& \cite{IAPWS1997}\\
      $\rho_h$	& \multicolumn{2}{c}{$900 $} & $kg\cdot m^{-3}$ 	& \cite{Sloan2007Book}\\
      $\rho_s$	& \multicolumn{2}{c}{$2100$} & $kg\cdot m^{-3}$ 	\\ 
% \\
% 
      \multicolumn{4}{c}{\textbf{Hydraulic properties}} \\
      $\lambda_{BC}$		& \multicolumn{2}{c}{$1.2$}		&		& \cite{Helmig1997}\\
      $P_{entry}$		& \multicolumn{2}{c}{$50$}		& $kPa$		& \cite{Helmig1997}\\ 
% \\
%      
      \multicolumn{4}{c}{\textbf{Hydrate kinetics}} \\
      \rowcolor{red!10}$k_{reac}$	& formation:	     &$0.2 \times 10^{-11}$			& mol$\cdot m^{-2}\cdot $	& \\
      \rowcolor{red!10}			& dissociation:	     &$3.2 \times 10^{-10}$			& $Pa^{-1}\cdot s^{-1}$			& \\
      $N_{Hyd}$			& \multicolumn{2}{c}{$5.75$}						& 					 \\
%       $P_{e,0}$	 		& \multicolumn{2}{c}{$ A_1 = 10^6$, $A_2 = 14.17$, $A_3 = 1886.79$}	& $Pa$ 				 	 \\ 38.98 8533.8
      $P_{e,0}$	 		& \multicolumn{2}{c}{$ A_1 = 10^6$, $A_2 = 38.98$, $A_3 = 8533.8$}	& $Pa$ 				 	& \cite{KamathHolder1987}\\
      $\dot Q_h$ 		& \multicolumn{2}{c}{$B_1 = 56599$, $B_2 = 16.744$}			& $W \cdot m^{-3}$ 			& \cite{Esmaeilzadeh2008} \\ 
% \\
% 
      \multicolumn{4}{c}{\textbf{Poroelasticity parameters}} \\
      $\alpha_{biot}$			& \multicolumn{2}{c}{$0.8$ }		& 	&\cite{Verruijt2008}\\
      $\nu_{sh}$			& \multicolumn{2}{c}{$0.15$}		& 	&\cite{Miyazaki2011}\\
      \rowcolor{red!10}			& \multicolumn{1}{c}{\cellcolor{red!10}formation}	& \multicolumn{1}{c}{\cellcolor{red!10}dissociation}		&   	 &\\
      \rowcolor{red!10}$E_{s}$		& \multicolumn{1}{c}{\cellcolor{red!10}$32$} 		& \multicolumn{1}{c}{\cellcolor{red!10}$160$}			& $MPa$	 &\\
      \rowcolor{red!10}$E_{h}$		& \multicolumn{1}{c}{\cellcolor{red!10}$250$}		& \multicolumn{1}{c}{\cellcolor{red!10}$360$}			& $MPa$	 &\\ 
      \rowcolor{red!10}$c$ 		& \multicolumn{1}{c}{\cellcolor{red!10}$1$}		& \multicolumn{1}{c}{\cellcolor{red!10}$3$}			& 	 &\\
\end{longtable}
% \twocolumn

\newpage

\paragraph{Acknowledgements}
We gratefully acknowledge the support for the first author by the German Research Foundation (DFG), through project no. WO 671/11-1. 
This work was further funded by the German Federal Ministries of Economy (BMWi) and Education and Research (BMBF) through the SUGAR project (grant No. 03SX250, 03SX320A \& 03G0856A), 
% by DEA Deutsche Erd\"ol AG, 
and the EU-FP7 project MIDAS (grant agreement no. 603418).
The reported experimental data is attached as supplemental material with this article. 

\bibliographystyle{plain}
\bibliography{myBib}

\end{document}